\magnification=\magstep1

\def\itemitemitem{\par\indent\quad
   \hangindent=3\parindent \textindent}

\null
\vskip 0.5 in
\centerline {\bf SNARKS FROM A K\'ASZONYI PERSPECTIVE: 
     A SURVEY}
     
\vskip 0.5 in

\noindent Richard C.\ Bradley \hfil\break
Department of Mathematics \hfil\break
Indiana University \hfil\break
Bloomington \hfil\break
Indiana 47405 \hfil\break
USA \hfil\break

\noindent bradleyr@indiana.edu \hfil\break

\vskip 0.5 in

   {\bf Abstract.}\ \ This is a survey or 
exposition of a particular collection of results 
and open problems involving snarks --- simple ``cubic'' 
(3-valent) graphs for which, for nontrivial reasons, 
the edges cannot be 3-colored.
The results and problems here are rooted in a series
of papers by L\'aszl\'o K\'aszonyi that were 
published in the early 1970s.
The problems posed in this survey paper can
be tackled without too much specialized 
mathematical preparation, and in particular seem well
suited for interested undergraduate mathematics 
students to pursue as independent research projects.
This survey paper is intended to facilitate research
on these problems. 
\hfil\break

\vfill\eject

\noindent {\bf 1.\ \ Introduction} \hfil\break

   The Four Color Problem was a major fascination for
mathematicians ever since it was posed in the middle
of the 1800s.  
Even after it was answered affirmatively (by 
Kenneth Appel and Wolfgang Haken in 1976-1977, with heavy
assistance from computers) and finally became transformed
into the Four Color Theorem, it has continued to be of
substantial mathematical interest.
In particular, there is the ongoing quest for a shorter
proof (in particular, one that can be rigorously checked  
by a human being without the assistance of a computer);
and there is the related ongoing philosophical 
controversy as to whether a proof that requires the use of 
a computer (i.e.\ is so complicated that it cannot be rigorously checked by a human being without the
assistance of a computer) is really a proof at all.
For a history of the Four Color Problem, a gentle 
description of the main ideas in the proof of the Four
Color Theorem, and a look at this controversy over whether
it is really a proof, the reader is referred to the book
by Wilson [Wi]. 
With the Four Color Problem having been ``solved'',
it is natural to examine related ``coloring''
problems that might have the kind of fascination
that the Four Color Problem has had.
\medskip

The Four Color Theorem has a well known equivalent
formulation in terms of colorings of edges of graphs:
If a graph $G$ has finitely many vertices and edges,
is ``cubic'' (``3-valent'', i.e.\ each vertex is 
connected to exactly three edges),
is planar, and satisfies certain other technical
conditions (to avoid trivial technicalities), then
it can be ``edge-3-colored'', i.e.\ its edges can each
be assigned one of three colors in such a way that no
two adjacent edges have the same color. \medskip
 
   Now if one removes the condition that $G$ be planar
but keeps all of the other conditions, then it may
be the case that $G$ {\it cannot\/} be edge-3-colored.
In a {\it Scientific American\/} article by the 
famous mathematics expositor Martin Gardner [Ga] in 1976,
the term ``snark'' was coined 
(or rather borrowed from Lewis Carroll's tale, 
{\it The Hunting of the Snark\/})
for such counterexamples $G$ for which the reason 
for the absence of an edge-3-coloring 
is in a certain technical sense ``nontrivial''.
(More on that in Definition 2.9 in 
Section 2 below.)
\medskip

   For a long time, only very few such ``nontrivial''
counterexamples had been known.
Then in 1975, Rufus Isaacs [Is] brought to light some 
infinite families of such nontrivial counterexamples, 
and some new methods for constructing such examples.
A year later, the 1976 paper of Gardner [Ga] alluded 
to above, popularized that work of Isaacs and
established the term ``snark'' for such examples. 
Since then, snarks have been a topic of extensive 
research by many mathematicians. \medskip

   However, there is a particular line of open questions 
on snarks which has received little attention during
that time.
Those open questions can be tackled
without too much specialized mathematical preparation.
In particular, those open questions seem to be reasonably
well suited for undergraduate mathematics students to 
pursue as independent research projects --- for example 
in an REU (Research Experience for Undergraduates) 
program in Mathematics.
This survey or expository paper will hopefully help 
anyone who becomes fascinated with snarks 
get started relatively quickly on research on this
particular collection of open problems. \medskip 

   The particular collection of research problems posed 
in this expository paper is rooted in three
obscure, relatively little known papers by 
L\'aszl\'o K\'aszonyi [K\'a1, K\'a2, K\'a3]
that were published in the early 1970s.
Further work on this collection of problems was done
later on by the author [Br1, Br2]; and 
still further work was done subsequently by 
Scott McKinney [McK] as part of a Mathematics REU program
when he was an undergraduate mathematics major. 
This expository paper is intended to enable the 
interested reader to quickly become acquainted with the relevant material in those six papers, and to get started quickly on the collection of research problems posed 
here.
\hfil\break  

   {\bf The organization of this paper.}\ \  
The collection of open problems is presented in
Section 7, the final section of this paper.
Sections 2 through 6 are intended as preparation
for those problems, and can be summarized as follows: 
 
Section 2 gives most of the relevant basic background 
material and terminology.
 
Section 3 gives a convenient, generously detailed,
(hopefully) easy-to-digest
presentation of the material in the three key 
papers of K\'aszonyi [K\'a1, K\'a2, K\'a3] that is
most directly relevant to the problems posed in Section 7.
Theorem 3.3, a result of K\'aszonyi [K\'a2, K\'a3],
is in essence the ``backbone'' of this entire survey
paper.
Theorem 3.5 deals with a very important classic
example --- the Petersen graph --- in connection with
Theorem 3.3.

Section 4 gives a detailed presentation of
some related material from [Br1, Br2] involving 
pentagons (cycles with five edges) in snarks.

Section 5 is intended to get the reader ``oriented'' 
to the collection of open problems posed in Section 7. 
In particular, it is
intended to motivate and --- through
Theorem 5.3 and its proof --- portray in a simple
context the use of certain techniques 
(in particular,  
Kochol [Ko1] ``superpositions'' for creating snarks)
that are pertinent to research on most of the problems 
posed in Section 7.  

Section 6 briefly explains the main results of
McKinney [McK] (involving techniques similar to
those in the proof of Theorem 5.3).
\hfil\break

   Much of this paper can be read quickly and somewhat superficially (without working through all of the 
detailed arguments).
However, to be able to develop the information and
techniques most important for effective work on the 
majority of problems posed in Section 7 (or on certain 
other, related ones), one needs to  
(1) master Section 2 (which is just background 
material and terminology and can be absorbed pretty 
quickly), \break 
(2) work through carefully, with pencil and paper, with 
the drawing of diagrams, the proofs of 
Theorems 3.3, 3.5, and 5.3, and
(3) read carefully all of the rest of Section 5 as well
(to become well oriented to the collection of problems
in Section 7).
\medskip

   Depending on what particular problems one chooses to 
work on (either ones posed in Section 7 or other, related ones), one might want to get hold of certain other 
references. 
In [Is], [Br2, Section 2], and [McK], one can find the 
proofs of certain results that are stated in this
survey paper without proofs. 
From the paper of Kochol [Ko1] itself, one can learn 
more about Kochol ``superpositions'', a key tool for 
some of the problems posed in Section 7.
For the Four Color Problem itself (from which sprouted
the study of snarks as well as much other mathematics),
a fascinating, easy-to-read historical account is given 
in the book by Wilson [Wi].  
\hfil\break

\noindent {\bf 2.\ \ Background material} \hfil\break

   All ``Remarks'' and informal comments, in this
section and throughout this paper, are well known,
standard, trivial facts, and their proofs will usually
be omitted.
\medskip

   Section 2 here is devoted to some background
terminology.
Most of it is standard, well established.
\medskip

   In this paper, a ``graph'' is always an undirected
graph.
It is always assumed to be nonempty 
(i.e.\ to have at least one vertex).
It is always assumed (typically with explicit reminder)
to have only finitely many vertices and have 
no loops and no multiple edges.
\medskip
 
Two or more graphs are (pairwise) ``disjoint'' 
(from each other) if no two of them have any vertex 
(or edge) in common.
\medskip

   Within a given graph, a vertex $v$ and an edge $e$ 
are said to be ``connected'' (to each other) if $v$ 
is an endpoint of $e$.
Two edges are said to be ``adjacent'' (to each other)
if they are connected to the same vertex.
Two vertices are said to be ``adjacent'' (to each other)
if they are connected to (i.e.\ are the endpoints of) 
the same edge.
\medskip

   Within a given graph, for a given positive 
integer $n$, a given vertex $v$ is said to be 
``$n$-valent'' if it is connected to exactly
$n$ edges.
A 1-valent vertex is also said to be ``univalent''.
\medskip

   The cardinality of a given set $S$ will be denoted
``${\rm card}\, S$\thinspace ''. 
\hfil\break

{\bf Definition 2.1.}\ \ Consider the group
${\bf Z}_2 \times {\bf Z}_2$ 
with the operation (denoted $+$) being coordinatewise 
addition modulo 2.
Let its elements be denoted as follows:
$$ 0 := (0,0), \quad a := (0,1), \quad b := (1,0), 
\quad c := (1,1)\ . \eqno (2.1) $$
Thus $0$ is the identity element, $a+a = b+b = c+c = 0$, 
also $c + b = a$, and so on. \hfil\break

   In this paper, in all ``edge-3-colorings'' 
(Definition 2.4(b) below), the ``colors'' will be the nonzero elements $a$, $b$, and $c$ of ${\bf Z}_2 \times {\bf Z}_2$ 
as defined in (2.1).
The group structure of ${\bf Z}_2 \times {\bf Z}_2$ itself
has for a long time been a handy ``bookkeeping'' tool in connection with edge-3-colorings 
(see e.g.\ its pervasive use in [Ko1])
and also in connection with the Four Color Theorem 
(see [Wi]).  
Its use will be illustrated in a few places below, and
especially in the proof of Theorem 5.3 in Section 5.
\medskip  

   It is tacitly understood that the use of those three
``colors'' $a$, $b$, and $c$ in (2.1) is only for
convenience, and that the results on edge-3-colorings 
given in this paper trivially carry over to the use of 
any other choice of three ``colors''.
\bigskip

   {\bf Remark 2.2.}\ \ If $x$, $y$, and $z$ are each an element of the set $\{a, b, c\}$ 
(the set of non-zero elements of 
${\bf Z}_2 \times {\bf Z}_2$ as defined in (2.1)), then the following two statements are equivalent:
\hfil\break
(i) $x$, $y$, and $z$ are distinct (that is, $\{x,y,z\} = \{a,b,c\}$);
\hfil\break
(ii) $x + y + z = 0$ (the element $(0,0)$ of 
${\bf Z}_2 \times {\bf Z}_2$). \hfil\break

{\bf Definition 2.3.}\ \ Suppose $G$ is a (not necessarily connected) graph that has only finitely many vertices, no loops, and no multiple edges. \medskip

   (a) The set of all edges of $G$ will be denoted 
${\cal E}(G)$. \medskip

   (b) The set of all vertices of $G$ will be denoted 
${\cal V}(G)$. \medskip

   (c) A ``cycle'' in $G$ is of course a subgraph with 
(for some integer $n \geq 3$) distinct vertices 
$v_1, v_2, \dots, v_n$ and edges $(v_i, v_{i+1})$, 
$i \in \{1, 2, \dots, n-1\}$, and $(v_n, v_1)$.
For a given integer $n \geq 3$, an ``$n$-cycle'' is a 
cycle with exactly $n$ vertices.  
A $5$-cycle will also be called a ``pentagon''.
\medskip

   (d) Two or more given cycles in $G$ are (pairwise) ``disjoint'' (from each other) if no two of them have
any vertex in common. \hfil\break

   {\bf Definition 2.4.}\ \ Suppose $G$ is a (not 
necessarily connected) graph that has only finitely many vertices, no loops, and no multiple edges, and each 
vertex of $G$ has valence 1, 2, or 3.  (The valences of 
the different vertices are not assumed to be be equal.) \medskip

   (a) A ``3-edge-decomposition'' of $G$ is a partition of 
${\cal E}(G)$ into three classes of edges such that no two adjacent edges of $G$ are in the same class.  The set of 
all 3-edge-decompositions of $G$ is denoted $ED(G)$.  
For a given $\delta \in ED(G)$ and a given pair of
edges $e_1$ and $e_2$ of $G$, the notation $e_1 \sim e_2$ means that $e_1$ and $e_2$ belong to the same class 
(of the three classes) in the decomposition $\delta$.
\medskip

   (b) An ``edge-3-coloring'' of $G$ is a function 
$\gamma : {\cal E}(G) \to \{a,b,c\}$ (where $a$, $b$, 
and $c$ are the nonzero elements of 
${\bf Z}_2 \times {\bf Z}_2$ defined in (2.1)) 
such that $\gamma(e_1) \neq \gamma(e_2)$ for every 
pair of adjacent edges $e_1$ and $e_2$ of $G$.  
The set of all edge-3-colorings of $G$ is denoted $EC(G)$.   If the set $EC(G)$ is nonempty, then one says simply
that the graph $G$ is ``edge-3-colorable'' or 
``can be edge-3-colored.''
\medskip

   An edge-3-coloring of $G$ is sometimes simply
called a ``coloring'' of $G$.
The phrases ``edge-3-colorable'' and 
``can be edge-3-colored'' are sometimes simply  
abbreviated ``colorable'' and ``can be colored''.
If $G$ cannot be (edge-3-)colored, it is often
referred to simply as being ``uncolorable'' or
``noncolorable''.
\medskip

   Edge-3-colorings, of cubic graphs 
(Definition 2.5(b) below), were first studied
by Peter Guthrie Tait in 1880 (see [Wi, Chapter 6]),
and are often called ``Tait colorings''. 
\bigskip

   {\bf Definition 2.5.}\ \ Suppose $G$ is a (not necessarily connected) graph that has only finitely many vertices, no loops, and no multiple edges. \medskip 

   (a) This graph $G$ is said to be ``quasi-cubic'' if every vertex of it is either 1-valent (univalent) or 3-valent.  
\medskip

   (b) The graph $G$ is said to be ``cubic'' if every
vertex of it is 3-valent. \medskip  

   Thus if $G$ is cubic, then it is quasi-cubic. 
\hfil\break

   If $G$ is a quasi-cubic graph in which no 
vertex is 3-valent (i.e.\ every vertex is univalent), 
then $G$ trivially has just one
edge or just a collection of isolated edges.
Such graphs will ordinarily be explicitly excluded
(as in the next Remark).
\bigskip 
 
   {\bf Remark 2.6.}\ \ Suppose $G$ is a (not 
necessarily connected) graph with finitely many 
vertices, no loops, and no multiple edges, and $G$ 
is quasi-cubic. 
Suppose further that $G$ has at least one 
3-valent vertex $v$. 
Let $e_1$, $e_2$, and $e_3$ denote the three edges
connected to $v$.
Then for every $\delta \in ED(G)$ (if one exists)
there exists exactly 
one coloring $\gamma \in EC(G)$ such that
(i) the equalities  
$$ \gamma(e_1) = a,\ \ \gamma(e_2) = b,\ \
\ \gamma(e_3) = c  \eqno (2.2)  $$
hold (again see (2.1)) and 
(ii) the edges in any of the three classes in the
decomposition $\delta$ are assigned the same color by 
$\gamma$.
In fact this induces a one-to-one correspondence between
$ED(G)$ and the set of all $\gamma \in EC(G)$ such
that (2.2) holds.
Hence 
$${\rm card}\thinspace ED(G)\ =\  
{\rm card}\thinspace \{\gamma \in EC(G): (2.2)\ 
{\rm holds}\}.  \eqno (2.3) $$
Since there are 6 permutations of the three colors
$a$, $b$, and $c$, it follows that each $\delta \in ED(G)$
induces exactly 6 colorings $\gamma \in EC(G)$ such that
the edges in any of the three classes in $\delta$ are
assigned the same color by $\gamma$.
\medskip 

   {\it To summarize\/}: \quad If $G$ is a (not necessarily connected)
graph with finitely many vertices, no loops, and no
multiple edges, and $G$ is quasi-cubic and has at least
one 3-valent vertex, then
$$ {\rm card}\thinspace EC(G)\ =\ 
6 \cdot {\rm card}\thinspace ED(G).  \eqno (2.4) $$
(Of course this holds trivially, with both sides 
being 0, in the case where no 3-edge-decomposition of 
$G$ exists, i.e.\ when $G$ cannot be colored.)
\hfil\break

   {\bf Definition 2.7.}\ \ A graph $G$ (not necessarily
quasi-cubic) is said to be ``simple'' if it has only 
finitely many vertices, has no loops and no multiple 
edges, and is {\it connected\/}. \hfil\break

   {\it Remark.}\ \ If a graph is both simple and
quasi-cubic and has at least one 3-valent vertex, then
trivially no edge of it can be connected to two
univalent vertices.
That is to be tacitly kept in mind in what follows.
\hfil\break

   {\bf Definition 2.8.}\ \ Suppose $G$ is a simple graph
(not necessarily quasi-cubic). \medskip

   (a) If $Q$ is a proper subset of ${\cal V}(G)$, 
then $G-Q$ denotes the graph that one obtains from 
$G$ by deleting every vertex in $Q$ and every edge 
that is connected to at least one vertex in $Q$. 
(All vertices in the set ${\cal V}(G) - Q$ are
left intact, even the ones that are endpoints of edges that
are removed.) \medskip 

   (b) If $S$ is any subset of ${\cal E}(G)$, then $G-S$ denotes the graph that one obtains from $G$ by deleting
every edge in $S$ (but not deleting any vertex).   
That is, $G-S$ is the graph that consists of 
(i) all vertices of $G$ and 
(ii) all edges of $G$ except the ones in $S$. \medskip

   (c)  If the graph $G$ has at least one cycle,
then the ``girth'' of $G$ is the smallest integer $n$ 
($\geq 3$) such that $G$ has an $n$-cycle.
\smallskip
 
(If $G$ has no cycles, i.e.\ $G$ is a tree, then its 
``girth'' is not defined.) \medskip 

   (d)  Suppose the graph $G$ has at least one pair of
disjoint cycles (recall Definition 2.3(d)).
For a given integer $n \geq 2$, the graph $G$ is 
said to be (at least) ``cyclically $n$-edge-connected'' if 
the following holds:  If $C_1$ and $C_2$ are any two 
disjoint cycles in $G$, and $S$ is any subset of 
${\cal E}(G)$ with ${\rm card}\ S \leq n-1$ such that $S$ 
has none of the edges in the cycles $C_1$ and $C_2$
(edges in $S$ may have {\it endpoint vertices\/} in $C_1$
and/or $C_2$), then there exists in the graph $G-S$ a 
path from (a vertex of) $C_1$ to (a vertex of) $C_2$.
\smallskip 

   (For a simple graph $G$ that does not have any pair of disjoint cycles, the notion of  
``cyclically $n$-edge-connected'' is not defined.) 
\medskip

   In graph theory, there are other notions of
``$n$-connectedness'' (for a given positive
integer $n$); but the one in (d) (``cyclically
$n$-edge-connected'') is the one that is most directly
relevant to the material in this survey, and is the one
that we shall stick with.  
\medskip

   (e) The following two facts are well known and
elementary to verify, and they are relevant to the next
definition.
(i) If a $G$ is a simple graph whose vertices all have
valence at least 2, then $G$ has at least one cycle.
(ii) If $G$ is a simple {\it cubic\/} graph whose girth is
at least 5, then it has at least one pair of disjoint
cycles --- in fact, for any cycle $C$ with the 
{\it minimum\/} 
number of vertices (namely the girth), there exists
another cycle that is disjoint from $C$.
(For (ii), under the conditions on $G$ and $C$, no
vertex of $G - {\cal V}(C)$ can be adjacent (in $G$) 
to more than
one vertex of $C$; and as a consequence, (i) implies
that $G - {\cal V}(C)$ has a cycle.)  
\hfil\break     

   {\bf Definition 2.9 (snarks).}\ \ Suppose $G$ is a simple cubic graph. 
Then $G$ is a ``snark'' if it has all of the following 
three properties: \hfil\break
(i) the girth of $G$ is at least $5$, \hfil\break
(ii) $G$ is (at least) cyclically 4-edge-connected, 
and \hfil\break
(iii) $G$ cannot be edge-3-colored. \medskip

   {\it Remark.}\ \ Here the ``primary'' property is (iii). 
The ``secondary'' properties (i) and (ii) are standard
restrictions, and in essence their purpose is to exclude
noncolorable simple cubic graphs 
whose noncolorability can be ``reduced''
in a trivial way to that of some
smaller noncolorable simple cubic graph.  
The term ``snark'' was suggested (borrowed from
Lewis Carroll) by Martin Gardner [Ga] as a convenient 
term to describe the class of graphs that had 
been studied in the paper published a year earlier by
Rufus Isaacs [Is] --- the graphs meeting the
conditions in Definition 2.9. 
\hfil\break

{\bf Remark 2.10 (the Petersen graph).}\ \ (i) The
classic ``smallest''
example of a snark is the Petersen graph, henceforth 
denoted ${\cal P}$, with ten vertices $u_i$ and $v_i$, 
$i \in \{0,1,2,3,4\}$ and fifteen edges $(u_i, u_{i+1})$, 
$(u_i,v_i)$, and $(v_i, v_{i+2})$, $i \in \{0,1,2,3,4\}$ 
(the addition in the subscripts is modulo 5).
\medskip

   (The proof that the Petersen graph ${\cal P}$ is a
snark, is straightforward.  
One can show that it has no 3-cycles or 4-cycles, and
that any two {\it disjoint\/} cycles must be two
pentagons (5-cycles) directly connected by the 
remaining 5 edges.
To see that ${\cal P}$ cannot be edge-3-colored, one
can first identify possible color patterns for the
edges $(u_i, u_{i+1})$, $i \in \{0,1,2,3,4\}$
--- using symmetries and permutations of colors, one can
reduce to basically just one color pattern --- and one
can then see that such a color pattern dictates the 
colors for the five edges $(u_i,v_i)$, and those
colors in turn lead to a contradiction when one tries
to assign colors to the remaining edges 
$(v_i, v_{i+2})$.)
\medskip 
        
   (ii) The Petersen graph ${\cal P}$ is loaded with symmetries.
Here is just a small part of that story:
If $P_1$ and $P_2$ are any two pentagons (5-cycles) 
in ${\cal P}$ (those two pentagons may be identical or
overlapping or disjoint), and one maps $P_1$ onto $P_2$
in any of the (ten) possible ways, then that mapping 
extends to a unique automorphism of ${\cal P}$. 
Also, every edge of ${\cal P}$ belongs to a pentagon
(in fact, to four pentagons).
It follows that any edge of ${\cal P}$ can be mapped to
any other edge as part of some automorphism of ${\cal P}$.
\medskip

   In papers such as those of Isaacs [Is] and 
Kochol [Ko1], methods have been devised to ``combine'' 
two or more snarks in order to form a ``bigger'' 
snark --- and thereby to create recursively, starting
with (say) the Petersen graph, infinite classes of 
arbitrarily large snarks.  
More on that later in this survey. \medskip

   Terms such as ``a Petersen graph'' (singular)
or ``Petersen graphs'' (plural) will be used for
referring to one or more graphs that are isomorphic 
to ${\cal P}$.
\hfil\break  

{\bf Definition 2.11.}\ \ Suppose $H$ is a simple 
quasi-cubic graph that has at least one 3-valent
vertex and can be edge-3-colored. \medskip

   (a) Suppose $\gamma \in EC(H)$.  Suppose $K$ is a 
subgraph of $H$ such that (i) $K$ is connected, 
(ii) $\gamma$ assigns just two colors, say $x$ and $y$, 
to the edges in $K$, and (iii) all edges of $H$ that are 
connected to (vertices of) $K$ but are themselves not 
in $K$, are assigned the third color, say $z$.  
Then $K$ is called a ``Kempe chain'' for the 
edge-3-coloring $\gamma$.
(That is, a Kempe chain is a subgraph $K$ which is
``maximal'' with respect to having 
properties (i) and (ii).) \medskip

   (b) In (a), if the two colors assigned to the edges in $K$ are $x$ and $y$, then $K$ is also called more specifically an $xy$-Kempe chain (for the coloring $\gamma$).
\smallskip 
   
   (Note that for a given coloring $\gamma$ and two
given colors $x$ and $y$, any two different $xy$-Kempe
chains must be disjoint from each other.
Of course an $xy$-Kempe chain can intersect an
$xz$-Kempe chain, where $z$ is the third color.)
\medskip

   (c) If a Kempe chain (for a given $\gamma \in EC(H)$) is a cycle, then it is called a ``Kempe cycle'' (or an 
``$xy$-Kempe cycle'' if the two colors are $x$ and $y$). \medskip

   (d) {\it Remark.}\ \ It is easy to see that for a given
Kempe chain $K$ (for a given $\gamma \in EC(H)$), either 
(i) $K$ is simply a path with two distinct endpoints, 
each of which is univalent, or \break 
(ii) $K$ is a (Kempe) cycle.  
If $H$ is cubic, then every Kempe chain must be a (Kempe) cycle.
\medskip

   (e) {\it Remark.}\ \ Suppose $\gamma$ is an 
edge-3-coloring of $H$, $x$ and $y$ are distinct colors
(in the set $\{a,b,c\}$ of colors from (2.1)), and 
$K$ is an $xy$-Kempe chain for $\gamma$.
  
(i) Suppose one interchanges the colors $x$ and $y$ on
just the Kempe chain $K$.
That is, suppose one defines the mapping 
$\mu : {\cal E}(H) \to \{a,b,c\}$ as follows:
$$ \mu(e)\ :=\ \cases {
x & if $e \in {\cal E}(K)$ and $\gamma(e) = y$ \cr
y & if $e \in {\cal E}(K)$ and $\gamma(e) = x$ \cr
\gamma(e) & for all edges $e$ of $H$ except the ones in $K$. \cr} $$
Then $\mu$ is an edge-3-coloring of $H$, and also $K$ is an
$xy$-Kempe chain for $\mu$ (as well as for $\gamma$).  

(ii) If one now starts with this new coloring $\mu$ and
one (again) interchanges the colors $x$ and $y$ on $K$,
then one obtains the original coloring $\gamma$.    
 \hfil\break

   {\bf Remark 2.12.}\ \ Here is a quick review of some 
well known further elementary facts.  
Suppose $G$ is a simple quasi-cubic graph that has
at least one 3-valent vertex. \medskip

   (a) Consider the ordered pairs $(v,e)$ where $v$ is a vertex of $G$ and $e$ is an edge of $G$ connected to $v$.
The number of such ordered pairs is even (since each edge generates two of them).  
Also, the total (even) number of such ordered pairs is the 
sum of the valences (1 or 3) of the vertices.   
It follows that {\it the number of vertices of $G$ is 
even\/} (since otherwise the sum of the valences of the
vertices would be odd, not even).  \medskip

   (b) Now suppose $\gamma \in EC(G)$, and $x$ is one of the colors ($a$, $b$, or $c$ --- see (2.1)).  
The number of vertices connected to an edge colored $x$ (by $\gamma$) is even (since each edge colored $x$
is connected to two such vertices).  
Hence the number of vertices {\it not\/}
connected to an edge colored $x$ is even, and of course 
such vertices must be univalent. \medskip

   (c) As a trivial consequence, for a given
$\gamma \in EC(G)$, the number of univalent vertices that 
{\it are\/} connected to an edge colored $x$ 
is going to be  

\noindent (i) even for each color $x$ ($\in \{a,b,c\}$)
if the total number of univalent vertices is even,

\noindent (ii) odd for each color $x$ ($\in \{a,b,c\}$)
if the total number of univalent vertices is odd.
\smallskip

\noindent That fact is known as the ``Parity Lemma''. 
\medskip

   (d) As a key special case of remark (c)(ii) above, 
if exactly five edges are connected to (different) 
univalent vertices, then for a given $\gamma \in EC(G)$,
one color is given to exactly three of those five edges,
and the other two colors are each given to exactly
one of those five edges.
\medskip    
 
   (e) Quasi-cubic graphs and the ``Parity Lemma'' have 
always played a ubiquitous role in the study of snarks.  
The ``Parity Lemma'' has the following (equivalent)
well known classic formulation in terms of the non-zero elements of ${\bf Z}_2 \times {\bf Z}_2$ (see (2.1)): 
\hfil\break

   {\bf Lemma 2.13 (Parity Lemma).} \quad {\sl Suppose 
$H$ is a simple quasi-cubic graph which has at least one
3-valent vertex but is not cubic.
Let $e_1, e_2, \dots, e_n$ denote the edges of $H$ that
are each connected to a univalent vertex.

   Suppose also that $H$ is edge-3-colorable, and 
$\gamma$ is an edge-3-coloring of $H$.  
Then
$$ \sum_{i=1}^n \gamma(e_i) = 0. \eqno (2.5) $$
(Here the addition is that of 
${\bf Z}_2 \times {\bf Z}_2$ as in Definition 2.1;
and the right hand side is of course the zero 
element (0,0) of ${\bf Z}_2 \times {\bf Z}_2$ as 
in (2.1).)} 
\hfil\break

   {\it Remark.}\ \ For a coloring $\gamma$ of a 
colorable simple cubic graph $G$, if
$\{e_1, e_2, \dots, e_n\}$ is (in the sense of inclusion)
a ``minimal cut set'' of edges of $G$, then (2.5) holds.
To see that, one can simply apply Lemma 2.13 itself to 
either one of the two disjoint simple quasi-cubic graphs 
that one obtains from $G$ by (say) ``cutting each of the 
edges $e_i$ in half'' (and inserting a vertex 
at the end of each of the resulting ``strands'').  
(Here of course the term ``cut set'' means that the 
graph $G - \{e_1, \dots, e_n\}$ is not connected.)
\hfil\break
       
   Kochol [Ko2] gives a generalization of Lemma 2.13
and uses it in the study of ``graph coloring'' problems
in more general contexts than just cubic graphs.
That will not be treated further here.   
\hfil\break

   {\bf Remark 2.14 (the flower snarks).}\ \ Here 
we just point out a special class of snarks known as the 
``flower snarks'':
\medskip

   Suppose $n$ is an {\it odd\/} integer such that
$n \geq 5$.
Let $J_n$ denote the simple cubic graph that consists 
of $4n$ vertices\ 
$t_k,\ u_k,\ v_k,\ w_k,\ k \in \{0,1,2,\dots,n-1\}$
and the following $6n$ edges:
$(t_k, t_{k+1}),\ (u_k, v_{k+1}),\ (v_k, u_{k+1}),\
(w_k, t_k),\ (w_k, u_k),\ (w_k, v_k),\ 
k \in \{0,1,2,\dots,n-1\}$. 
Here addition in the subscripts is modulo $n$;
for example, for $k = n-1$, $t_{k+1} := t_0$.
\medskip

   Isaacs [Is, Theorem 4.1.1] showed that for each odd 
integer $n \geq 5$, this graph $J_n$ --- as Isaacs [Is] himself called it --- is a snark.
(We shall not give the argument here.)\ \
This class $J_n,\ n \in \{5, 7, 9, 11, \dots\}$
was one of the large classes of snarks that was
brought to light by Isaacs [Is].
These snarks $J_n,\ n \in \{5, 7, 9, 11, \dots\}$
were later referred to by Gardner [Ga] as the 
``flower snarks'', and that name has become customary 
since then for these particular snarks.
In what follows, the flower snarks will occasionally
be alluded to briefly.
It will just be mentioned here in passing that the 
flower snarks have special properties and a special fascination and are a topic of interest in 
their own right.
For example, Tinsley and Watkins [TW] computed the 
(topological) genus of each of the flower snarks.
\medskip

   Isaacs [Is] (and Gardner [Ga]) included $n=3$ in the definition above.
Technically, $J_3$ is not a snark, because it has a 
``triangle'' (a 3-cycle); but (as Isaacs [Is] noted) if 
one ``contracts'' that triangle in $J_3$ to a single 
vertex, one obtains the Petersen graph.
\medskip

   (To match the definition of the graphs $J_n$ given 
above to the pictures of those graphs in [Is, p.\ 233] 
and [Ga, p.\ 128], first decipher the ``Remark'' (both
paragraphs of it) in [Is, p.\ 234].)   
\hfil\break
\vfill\eject

\noindent {\bf 3.\ \ Results of K\'aszonyi} \hfil\break

   Throughout the rest of this survey paper, it is 
tacitly understood that theorems or proofs or informal comments that involve ``snarks'', 
often trivially carry over to --- and in
some cases were originally formulated for --- some
broader classes of noncolorable simple cubic graphs
that were allowed to satisfy less stringent
``secondary'' conditions than conditions (i) and (ii)
in Definition 2.9.
As a convenient formality, our discussion of noncolorable 
simple cubic graphs will be confined to snarks
(i.e.\ satisfying all three conditions in Definition 2.9); 
we shall thereby avoid having to occasionally bother with 
some trivial extra technicalities.
\hfil\break
 
   Section 3 here is an exposition of some of the work of 
K\'aszonyi [K\'a1, K\'a2, K\'a3] that pertains to 
snarks and also involves a closely related topic:
edge-3-colorable simple cubic graphs with 
``orthogonal edges'' (Definition 3.1 below).
The connection between those two topics is given as
part of Theorem 3.3 below.
In its entirety, Theorem 3.3 itself, a result of 
K\'aszonyi [K\'a2, K\'a3], is the ``backbone'' of this 
entire survey paper. 
\hfil\break

   {\bf Definition 3.1} 
(K\'aszonyi [K\'a1, K\'a2, K\'a3]).\ \  
Suppose $H$ is a simple cubic graph which can be
edge-3-colored.  
Suppose $d_1$ and $d_2$ are edges of $H$. 
Those edges $d_1$ and $d_2$ are said to be ``orthogonal'' 
(in $H$) if there does not exist $\gamma \in EC(H)$
for which $d_1$ and $d_2$ are edges in the same 
(two-color) Kempe cycle. \medskip

   {\it Remark\/.}\ \ Of course in the context of this definition, the
``orthogonal'' edges $d_1$ and $d_2$ cannot be adjacent.
(If $d_1$ and $d_2$ were adjacent, then for a given 
edge-3-coloring 
$\gamma$ of $H$, if one lets $x$ and $y$ denote the colors of $d_1$ and $d_2$ respectively, those two edges would 
belong to the same
$xy$-Kempe cycle for $\gamma$, contradicting the 
definition of ``orthogonal edges''.) \hfil\break

{\bf Notations 3.2.}\ \ Suppose $G$ is a simple cubic graph that has girth at least 4 and is at least cyclically 
2-edge-connected, and $e = (u,v)$ is an edge of $G$.   Let $t_i,\ i \in \{1,2\}$ denote the two vertices $\neq v$ 
that are adjacent to $u$ in $G$; and
let $w_i,\ i \in \{1,2\}$ denote the two vertices $\neq u$ that are adjacent to $v$ in $G$.  
(By the assumptions on $G$, those four vertices $t_1$, $t_2$, $w_1$, and $w_2$ are distinct, and $G$ does
not have an edge of either the form $(t_1,t_2)$ or 
$(w_1,w_2)$.)

   Let $G_e$ denote the simple cubic graph that consists
of the graph $G-\{u,v\}$ and the two new edges
$d_1 := (t_1,t_2)$ and $d_2 := (w_1,w_2)$.

  That is, to obtain $G_e$, one deletes from $G$ the two vertices $u$ and $v$ and all five edges (including $e$) connected to them, and one then inserts the two new edges $d_1$ and $d_2$.  
The simple cubic graph $G_e$ is conformal to the graph 
$G-\{e\}$. 
(The fact that $G_e$ is simple and cubic, is an elementary consequence of the assumptions here on $G$ itself.)  
\hfil\break

   In Notations 3.2, the graph $G$ is not assumed to be a snark.
However, apart from Problem 9 in Section 7.1, all of the applications of Notations 3.2 in this survey paper will explicitly involve the case where $G$ is 
assumed to be a snark. 
\hfil\break
\vfill\eject

   {\bf Theorem 3.3} (K\'aszonyi).\ \ {\sl (A) Suppose $H$
is a simple cubic graph which can be edge-3-colored, 
and $d_1$ and $d_2$ are orthogonal edges of $H$ 
(see Definition 3.1).
Then there exists a positive integer
$J$ such that the following three statements (i), (ii),
(iii) hold:
\hfil\break
(i) ${\rm card}\thinspace ED(H) = 3J$ (and hence
${\rm card}\thinspace EC(H) = 18J$). \hfil\break
(ii) ${\rm card}\thinspace
\{\delta \in ED(H): d_1 \sim d_2 \}\ =\ J$. \hfil\break
(iii) If $x$ and $y$ each $\in \{a,b,c\}$ ($x$ and
$y$ may be the same or different), then \hfil\break
${\rm card}\thinspace \{\gamma \in EC(H): \gamma(d_1) = x\
{\rm and}\ \gamma(d_2) = y \}\ =\ 2J$.

   (B) Suppose $G$ is a snark and $e$ is an edge of $G$.
Let $d_1$ and $d_2$ be the edges of $G_e$ specified in
Notations 3.2.
If the (simple cubic) graph $G_e$ is edge-3-colorable,
then $d_1$ and $d_2$ are orthogonal edges of $G_e$.

   (C) Suppose $G$ is a snark and $e$ is an edge of $G$.
Let $d_1$ and $d_2$ be the edges of $G_e$ specified in
Notations 3.2.
Then there exists a nonnegative integer $L$ such that the following three statements (1), (2), (3) hold: \hfil\break
(1) ${\rm card}\thinspace ED(G_e) = 3L$ (and hence
${\rm card}\thinspace EC(G_e) = 18L$). \hfil\break
(2) ${\rm card}\thinspace \{\delta \in ED(G_e): d_1 \sim d_2 \}\ =\ L$. \hfil\break 
(3) If $x$ and $y$ each $\in \{a,b,c\}$ ($x$ and $y$ may be the same or
different), then \hfil\break 
${\rm card}\thinspace \{\gamma \in EC(G_e): \gamma(d_1) = x\ {\rm and}\ \gamma(d_2) = y \}\ =\ 2L$.} \hfil\break

   This theorem, and its proof given below, are due to
K\'aszonyi [K\'a2, K\'a3].
In the proof of statement (A) given below, there are three
special classes of edge-3-colorings of $H$ that are called
(in the proof below) $Q_a$, $Q_b$, $Q_c$ 
(where $a,b,c$ are the ``colors'' from (2.1)).
The one-to-one correspondences (involving just the
interchange of the colors on particular Kempe cycles)
between those three classes, as spelled out in the proof
below, were pointed out by K\'aszonyi [K\'a3,
p.\ 35, the paragraph after Theorem 4.3
(together with p.\ 28, Section 1.15)] (in a somewhat
cryptic manner --- see Section 7.2(D) in Section 7).
Those one-to-one correspondences together give 
statement (A).
The proof of statement (B) given below is
taken from [K\'a2, p.\ 125, lines 8-24].
(K\'aszonyi's own formulation and argument there were 
technically slightly
more general, in the sense that the simple cubic
graph $G$, while assumed to be noncolorable, was allowed
to satisfy somewhat less stringent ``secondary'' 
conditions than conditions (i) and (ii) in 
Definition 2.9.
The term ``snark'' as in Definition 2.9 had not yet
been coined or codified.)\ \
Statement (C) is (aside from the trivial case where $G_e$
cannot be edge-3-colored) simply an application of
statements (A) and (B).
\hfil\break

   {\bf Proof.}\ \ {\bf Proof of statement (A).}\ \ Let 
$v$ be one
of the end-point vertices of the edge $d_1$, and let
$e_1$ and $e_2$ be the other two edges (besides $d_1$)
that are connected to $v$.
By the Remark after Definition 3.1, neither $e_1$ nor $e_2$
is the edge $d_2$.

   As described in Remark 2.6, every $\delta \in ED(H)$ induces a unique $\gamma \in EC(H)$ such that
$$ \gamma(d_1) = a,\ \ \gamma(e_1) = b,\ \ {\rm and}\ \ 
\gamma(e_2) = c.  \eqno (3.1) $$
For each $x \in \{a,b,c\}$, let $Q_x$ denote the set of
all $\gamma \in EC(H)$ such that (3.1) holds and
$\gamma(d_2) = x$.

   At this point, we make no assumptions on whether or not
the set $Q_x$, for a given color $x$, may be empty.
Accordingly, for now, definitions and arguments below
involving $Q_x$ are allowed to be ``vacuous''. 

   The assumption (in statement (A)) that the edges
$d_1$ and $d_2$ are orthogonal, has the following
elementary consequence:
For any $\gamma$ in either $Q_a$ or $Q_b$ (if such a 
$\gamma$ exists), if one interchanges the colors $a$ and
$b$ on the $ab$-Kempe cycle containing the edge $d_2$,
that will not change the colors (in (3.1)) of any of the 
edges $d_1$, $e_1$, or $e_2$.
Thus one has (for now, possibly vacuous) mappings 
$M: Q_a \to Q_b$ and $M^*: Q_b \to Q_a$ defined as
follows:
For each $\gamma \in Q_a$ (resp.\ $\gamma \in Q_b$), let
$M\gamma$ (resp.\ $M^*\gamma$) denote the element of
$Q_b$ (resp.\ of $Q_a$) which is obtained from $\gamma$
by the interchanging of the colors $a$ and $b$ on the 
$ab$-Kempe cycle (for $\gamma$) containing the edge $d_2$.
By a trivial argument (recall Remark (e)(i)(ii) in 
Definition 2.11), the mappings $M$ and $M^*$ are
inverses of each other:  $M^*M\gamma = \gamma$ for
$\gamma \in Q_a$, and $MM^*\gamma = \gamma$ for
$\gamma \in Q_b$.
It follows that $M$ is one-to-one and onto, as a mapping
from $Q_a$ to $Q_b$.
(And of course an analogous comment holds for $M^*$.)
Hence 
${\rm card}\thinspace Q_a = {\rm card}\thinspace Q_b$.

   By an exactly analogous argument,
${\rm card}\thinspace Q_a = {\rm card}\thinspace Q_c$. 
Thus
$$ {\rm card}\thinspace Q_a\ =\      
{\rm card}\thinspace Q_b\ =\ 
{\rm card}\thinspace Q_c.  \eqno (3.2) $$

   Define the nonnegative integer $J$ by
$$ J\ :=\ {\rm card}\thinspace Q_a.  \eqno (3.3) $$

   In the same manner as in Remark 2.6,
one has a one-to-one correspondence between the sets
$\{\delta \in ED(H): d_1 \sim d_2\}$ and $Q_a$.
(Every $\delta \in ED(H)$ such that $d_1 \sim d_2$,
trivially induces a unique $\gamma \in Q_a$.)\ \ 
Hence by (3.3), the equality in sub-statement (ii) in statement (A) holds.

   Also, precisely as in eq.\ (2.3) in Remark 2.6,  
$${\rm card}\thinspace ED(H)\ =\  
{\rm card}\thinspace \{\gamma \in EC(H): (3.1)\ 
{\rm holds}\}.  \eqno (3.4) $$
Now the set in the right hand side of (3.4) is simply
$Q_a \cup Q_b \cup Q_c$. 
Also, trivially by definition, those sets $Q_a$, $Q_b$,
and $Q_c$ are (pairwise) disjoint.  
Hence by (3.2), (3.3), and (3.4),
${\rm card}\thinspace ED(H) = 3J$.  Hence by
eq.\ (2.4) in Remark 2.6, 
sub-statement (i) in statement (A) holds.
Also, by the hypothesis of statement (A), the set
$EC(H)$ (or $ED(H)$) is nonempty.  
Hence by sub-statement (i) itself in statement (A), the
integer $J$ is positive.

   Hence, in the proof of statement (A), all that now 
remains is to verify sub-statement (iii), 
which is just an automatic, trivial by-product of 
sub-statements (i) and (ii).
Here are the details: 
\medskip

   {\it Proof of sub-statement (iii).}\ \ First, any 
$\delta \in ED(H)$ such that $d_1 \sim d_2$, 
induces six colorings $\gamma \in EC(H)$ such that
$\gamma(d_1) = \gamma(d_2)$ (with the edges in any of the
three classes in $\delta$ being assigned the same color
by $\gamma$).  
Hence by sub-statement (ii),
$$ {\rm card}\thinspace \{\gamma \in EC(H):
\gamma(d_1) = \gamma(d_2)\}\ =\ 6J.  \eqno (3.5) $$
By trivial permutations of colors, one can show that
for the three colors $x \in \{a,b,c\}$, the numbers
${\rm card}\thinspace \{\gamma \in EC(H): 
\gamma(d_1) = \gamma(d_2) = x \}$ 
are equal, and hence by (3.5) they must each be equal
to $2J$.
Thus sub-statement (iii) holds in the case where
$x = y$.

   Next, by (3.5) and sub-statement (i) (in statement (A)), 
$$ {\rm card}\thinspace \{\gamma \in EC(H):
\gamma(d_1) \neq \gamma(d_2)\}\ =\ 12J.  \eqno (3.6) $$   
By trivial permutations of colors, one can show that
for the six permutations $(x,y)$ of two distinct
elements of $\{a,b,c\}$, the numbers
${\rm card}\thinspace \{\gamma \in EC(H): \gamma(d_1) = x\ {\rm and}\ \gamma(d_2) = y \}$ 
are equal, and hence by (3.6) they must each be equal
to $2J$.
Thus sub-statement (iii) holds for the case $x \neq y$.
That completes the proof of sub-statement (iii), and of
statement (A). \medskip

   {\bf Proof of statement (B).}\ \ Suppose $G_e$ can
be edge-3-colored. 
Suppose the edges $d_1$ and $d_2$ are {\it not\/} orthogonal.
We shall aim for a contradiction.
The argument will be given here in a somewhat informal way, but that should not cause any confusion.
We shall make detailed use of the symbols in
Notations 3.2.

   Let $\gamma$ be an edge-3-coloring of $G_e$ such that
for two given colors, say $a$ and $b$, the edges $d_1$
and $d_2$ belong to the same $ab$-Kempe chain $K$.

   Let us ``reinsert'' the vertices $u$ and $v$ (see
Notations 3.2) into the middle of the edges $d_1$ and
$d_2$ respectively, thereby ``re-obtaining'' the graph
$G - \{e\}$ (which has two ``2-valent'' vertices
$u$ and $v$, and whose other vertices are still 3-valent).

Those two ``newly reinserted'' vertices $u$ and $v$
split the above-mentioned $ab$-Kempe chain $K$ 
for the coloring $\gamma$ into two pieces.  
On either one of those two pieces (but not the other), interchange the colors $a$ and $b$.
Thereby one obtains an edge-3-coloring of the graph
$G - \{e\}$ in which  
(i) no two adjacent edges have
the same color and (ii) the two edges connected to the
vertex $u$ (respectively to $v$) have (in either order)
the colors $a$ and $b$.

Now ``reinsert'' the edge $e$ to ``re-obtain'' the
original snark $G$, and assign to the ``newly
re-inserted'' edge $e$ the color $c$. 
One thereby obtains an edge-3-coloring of the snark $G$,
contradicting the definition of ``snark''.

Thus a contradiction has occurred, and hence the edges
$d_1$ and $d_2$ must be orthogonal after all.
That completes the proof of statement (B). \medskip

   {\bf Proof of statement (C).}\ \ In the case
where the graph $G_e$ itself cannot be edge-3-colored,
statement (C) (all parts of it) hold trivially with
$L = 0$.
In the case where the graph $G_e$ can be edge-3-colored,
statement (C) follows immediately from statements (A)
and (B) (with the graph $H$ in statement (A) being $G_e$
in statement (C), and with the integer $L$ in statement (C) being the integer $J$ in statement (A)).
That completes the proof of statement (C), and of 
Theorem 3.3.
\hfil ////\break

{\bf Definition 3.4.}\ \ Refer to Theorem 3.3(C)(1).  
For any snark $G$ and any edge $e$ of $G$, define the nonnegative integer $\psi(G,e)$ --- the 
``K\'aszonyi number'' of $G$ and $e$ --- as follows: 
$$\psi(G,e) := (1/3) \cdot {\rm card}\thinspace ED(G_e)\ .  \eqno (3.7) $$ 
That is, $\psi(G,e)$ denotes the nonnegative integer $L$
such that ${\rm card}\thinspace ED(G_e) = 3L$
(i.e.\ such that ${\rm card}\thinspace EC(G_e) = 18L$).
\hfil\break

   {\it Remark.}\ \ If $G$ is a snark, and the graph $H$
is an isomorphic copy of $G$, then $H$ is also a snark;
and if also $e$ is any edge of $G$, and $e^*$ is the edge
of $H$ which is the ``image'' of $e$ under a
given isomorphism from $G$ to $H$, then $\psi(H,e^*) = \psi(G,e)$.  The argument is elementary.
\hfil\break

   The notation $\psi(G,e)$, used in [Br2],
implicitly comes from the result of
K\'aszonyi [K\'a2, K\'a3] given in Theorem 3.3(C) above.
It will be used throughout the rest of this paper.
\hfil\break

   {\bf Theorem 3.5} (K\'aszonyi).\ \ {\sl For the Petersen 
graph ${\cal P}$ (see Remark 2.10) and any edge $e$ of
${\cal P}$, one has that $\psi({\cal P},e) = 1$.} 
\hfil\break

   This theorem is due to K\'aszonyi [K\'a1, K\'a2].
K\'aszonyi [K\'a2, p.\ 126, the Remark] pointed out
(in a slightly cryptic manner) that for the Petersen graph
${\cal P}$ and an edge $e$ of ${\cal P}$, the graph 
${\cal P}_e$ is (isomorphic to) a certain
graph consisting of an 
``8-vertex wheel with four `rim-to-rim spokes
going through the hub'\thinspace '' 
(with no vertex at the ``hub''). 
Slightly earlier, K\'aszonyi [K\'a1, pp.\ 81-82] had 
already shown that for that particular graph
(the ``8-vertex wheel with rim-to-rim spokes''), 
there are exactly three 3-edge-decompositions.  
Thereby (see (3.7)) Theorem 3.5 was established.
\medskip

   The proof given below is a slightly shortened version
of K\'aszonyi's argument, using the
terminology in Definition 3.4 and taking advantage
of statement (C)(3) in Theorem 3.3.
\hfil\break

   {\bf Proof.}\ \ By Remark 2.10(ii) and the Remark after
Definition 3.4, it suffices to carry out the argument for
any one particular edge $e$ of the Petersen graph $
{\cal P}$.
Referring to the formulation of the Petersen graph
${\cal P}$ in Remark 2.10(i), we shall consider the edge
$e := (u_2, u_3)$ there, and we shall express the graph
$W := {\cal P}_e$ as an 
``8-vertex wheel with rim-to-rim spokes''
(K\'aszonyi's [K\'a1, K\'a2] way of portraying 
${\cal P}_e$) as follows:
The eight vertices of $W = {\cal P}_e$ will be relabeled as
$t_i, i \in \{0, 1, 2, \dots, 7\}$, with
(see the notations in Remark 2.10(i))
$t_0:=u_4$, $t_1:=u_0$, $t_2:=u_1$, $t_3:=v_1$,
$t_4:=v_3$, $t_5:=v_0$, $t_6:=v_2$, and $t_7:=v_4$.
Then the twelve edges of $W$ fall into two
classes: the eight ``wheel edges'' 
$\epsilon_i := (t_i, t_{i+1})$, $i \in \{0, 1, \dots, 7\}$
(with addition mod 8 --- thus $\epsilon_7 := (t_7, t_0)$);
and the four ``spoke edges'' 
$f_i := (t_i, t_{i + 4})$, $i \in \{0,1,2,3\}$.
In the terminology of Notations 3.2, the edges $d_1$ and
$d_2$ are (say) respectively
$$\eqalignno{ d_1 &:= (u_4, v_3) = (t_0, t_4) = f_0 \quad
{\rm and} \cr
d_2 &:= (u_1, v_2) = (t_2, t_6) = f_2. & (3.8) \cr} $$ 

   Referring to (3.7), our task is to show that
${\rm card}\thinspace ED(W) = 3$.
Referring to (3.8) and Theorem 3.3(C)(1)(3), one has
that it suffices to show that
$$ {\rm card}\thinspace \{\gamma \in EC(W):
\gamma(f_0) = \gamma(f_2) = a\}\ =\ 2.  \eqno (3.9) $$  

   Let us count the ways of constructing colorings
$\gamma$ of $W$ such that $\gamma(f_0) = \gamma(f_2) = a$.
Each of the eight ``wheel edges'' $\epsilon_i$,
$i \in \{0, 1, \dots, 7\}$ is connected to one of the
vertices $t_0$, $t_2$, $t_4$, or $t_6$, the 
endpoints of the edges $f_0$ and $f_2$.
Hence those eight ``wheel edges'' must be colored
alternately $b$ or $c$.
There are two ways of doing that, depending on (say)
which color ($b$ or $c$) is assigned to the edge
$\epsilon_0$.
Either way, that forces one to complete the coloring
$\gamma$ of $W$ by assigning the color $a$ to the
remaining two ``spoke edges'' $f_1$ and $f_3$.
Thus (3.9) holds.
That completes the proof of Theorem 3.5.   
\hfil ////\break

   The rest of the material here in Section 3 is somewhat peripheral to the main theme of this survey paper, and is included here only to round out the picture a little.
Theorem 3.8 below is of independent interest.
\hfil\break

{\bf Notations 3.6.}\ \ Suppose $H$ is a simple cubic graph 
which can be edge-3-colored, and $d_1$ and $d_2$ are orthogonal
edges of $H$ (see Definition 3.1).  Let 
${\cal Q}(H, d_1, d_2)$
denote the set of all subgraphs $\Lambda$ of $H$ with the
following four properties: \hfil\break 
(i) $\Lambda$ is the union of one or more (pairwise) disjoint cycles in $H$.
\hfil\break
(ii) Each cycle in $\Lambda$ has an even number of edges. \hfil\break
(iii) $\Lambda$ contains every vertex of $H$. \hfil\break
(iv) Neither $d_1$ nor $d_2$ is an edge of $\Lambda$. 

   For each $\Lambda \in {\cal Q}(H, d_1, d_2)$, 
let ${\cal N}(\Lambda)$ denote the number of cycles in 
$\Lambda$. \hfil\break

{\bf Theorem 3.7} (K\'aszonyi).\ \ {\sl (A) Suppose $H$ is a simple cubic graph 
which can be edge-3-colored, and $d_1$ and $d_2$ are orthogonal
edges of $H$ (see Definition 3.1).  Then (see Notations 3.6 above)
$$ {\rm card}\thinspace ED(H)\ =\ (3/2) \cdot
\sum_{\Lambda \in {\cal Q}(H, d(1), d(2))} 2^{{\cal N}(\Lambda)}  
\eqno (3.10)$$     
(where $d_1$ and $d_2$ are written here as 
$d(1)$ and $d(2)$ for typographical convenience).

   (B) If $G$ is a snark, $e$ is an edge of $G$, and the (simple cubic) graph $G_e$ has no Hamiltonian cycle, then 
the nonnegative integer $\psi(G,e)$ is even.} \hfil\break

   This theorem, and its proof given below, are due to
K\'aszonyi [K\'a3, p.\ 35, Section 4.2].
Statement (A) was shown there --- together with some
other closely related information --- in an
indirect, somewhat cryptic form (see Section 7.2(D) 
in Section 7).
Statement (B) is simply a special case of statement (A).
\hfil\break

   {\bf Proof.}\ \ {\bf Proof of statement (A).}\ \ 
Refer to Theorem 3.3(A)(i)(iii) (and to Definition 2.1). 
To prove (3.10) and thereby statement (A),  
it suffices to show that
$$ {\rm card}\thinspace
\{\gamma \in EC(H): \gamma(d_1) = \gamma(d_2) = a\}\ =\ 
\sum_{\Lambda \in {\cal Q}(H, d(1), d(2))} 2^{{\cal N}(\Lambda)}\ . \eqno (3.11)$$

   For any $\gamma \in EC(H)$ such that 
$\gamma(d_1) = \gamma(d_2) = a$, the set of edges that
are colored $b$ or $c$ (by $\gamma$), together with
their endpoints, form a subgraph $\Lambda$ in the class
${\cal Q}(H, d_1, d_2)$.
Hence to count the colorings $\gamma \in EC(H)$ such
that $\gamma(d_1)=\gamma(d_2) = a$,  it suffices to take
the sum, over all $\Lambda \in {\cal Q}(H,d_1,d_2)$, of
the number of colorings $\gamma$ (with
$\gamma(d_1)=\gamma(d_2) = a$) whose $bc$-Kempe cycles
together form the graph $\Lambda$.
For each such $\Lambda$, there are
exactly $2^{{\cal N}(\Lambda)}$
such colorings, since for each cycle in $\Lambda$, there
are exactly two ways of assigning the alternating colors
$b$ and $c$.  Thus (3.11) holds.
That completes the proof of statement (A).
\medskip

   {\bf Proof of statement (B).}\ \ In the case where
the graph $G_e$ itself cannot be edge-3-colored, 
statement (B) is trivial, with $\psi(G,e) = 0$.
Therefore, assume $G_e$ can be colored.
 
   Let $H := G_e$, and 
let the edges $d_1$ and $d_2$ be as in Notations 3.2.
By Theorem 3.3(B), those two edges $d_1$ and $d_2$ are orthogonal.
If $\Lambda \in {\cal Q}(H, d_1, d_2)$, then it must have
at least two cycles (since otherwise $\Gamma$ itself would 
be a Hamiltonian cycle, contradicting the hypothesis
of statement (B)) --- that 
is ${\cal N}(\Lambda) \geq 2$ ---
and hence $2^{{\cal N}(\Lambda)}$ is a multiple of 4.
Hence the right side of (3.10) is a multiple of 6.  
Hence by (3.10) itself and (3.7), the integer
$\psi(G,e) = (1/3)\thinspace {\rm card}\thinspace ED(H)$ is even.
That completes the proof of statement (B), and of
Theorem 3.7.
\hfil ////\break

   {\bf Theorem 3.8} (K\'aszonyi).\ \ {\sl Suppose $H$ 
is a simple cubic graph which can be edge-3-colored, 
and it has at least one pair of orthogonal edges.  
Suppose further that \break 
${\rm card}\thinspace ED(H) = 3$ 
(the smallest possible number under the assumptions 
here --- see Theorem 3.3(A)).  
Then $H$ is nonplanar.} \hfil\break

   This theorem is due to K\'aszonyi [K\'a2, Theorem 6].  
This theorem and its proof are of intrinsic interest in 
their own right, but will not be needed anywhere else in 
this survey paper.
The proof is somewhat long and complicated and will not be 
repeated here. \hfil\break

\noindent {\bf 4.\ \ Pentagons} \hfil\break

   Recall from Definition 2.3(c) that a 5-cycle is also
called a ``pentagon''.
Not all snarks have a pentagon.
For example, the flower snarks 
$J_n,\ n \in \{7,9,11,\dots\}$ in Remark 2.14 do not
have a pentagon.
Kochol [Ko1] constructed snarks with arbitrarily large
girth. 
In contrast, the Petersen graph (see Remark 2.10) has 
twelve pentagons; and the flower snark $J_5$ 
(again see Remark 2.14) has exactly one.
Section 4 here will be narrowly focused on certain facts
involving pentagons in snarks.
\medskip
 
In particular, if $G$ is a snark and $P$ is a pentagon 
(if one exists) in $G$, then the numbers
$\psi(G,e),\ e \in {\cal E}(P)$ are equal. 
That and some related information will be given in 
Theorems 4.5 and 4.8 below, after some background 
information.  
Remarks 4.2 and 4.3 below are well known and quite 
trivial, and their proofs will be omitted. 
\hfil\break

{\bf Convention 4.1.}\ \ Here in Section 4, in notations 
such as $v_k$ for vertices, the indices $k$ will be taken 
as elements of the field ${\bf Z}_5$.  
The elements of ${\bf Z}_5$ will be denoted simply 
$0, 1, 2, 3, 4$, with addition and multiplication mod 5.

   For example, in such a context one has for $k = 4$ that 
$k+1 = 0$ and $2k = 3$. \hfil\break

{\bf Remark 4.2.}\ \ Suppose $x$, $y$, and $z$ are (in any ``order'') three distinct elements of ${\bf Z}_5$. \medskip 

   Consider the following three conditions:

\noindent (i) There exists $k \in {\bf Z}_5$ such that 
$\{x,y,z\} = \{k-1,k,k+1\}$. \hfil\break
(ii) There exists $k \in {\bf Z}_5$ such that 
$\{x,y,z\} = \{k-2,k,k+2\}$. \hfil\break
(ii$'$) There exists $j \in {\bf Z}_5$ such that
$\{2x, 2y, 2z\} = \{j-1, j, j+1\}$. \medskip

   Then by trivial arithmetic (mod 5), the following three statements hold: 
   
\noindent (A) Exactly one of conditions (i), (ii) holds. 
\hfil\break
(B) Conditions (ii) and (ii$'$) are equivalent. \hfil\break
(C) Hence, exactly one of conditions (i), (ii$'$) holds.
 \hfil\break

{\bf Remark 4.3.}\ \ Suppose $G$ is a simple cubic graph 
with girth 5,
and $P$ is a pentagon in $G$ with vertices $v_i$ and edges 
$(v_i, v_{i+1})$, $i \in {\bf Z}_5$.    
The graph $G - {\cal E}(P)$ is quasi-cubic, and its 
univalent vertices are precisely the ones 
$v_i,\ i \in {\bf Z}_5$.
For each $i \in {\bf Z}_5$, let $u_i$ denote the 
vertex of $G$ such that $\epsilon_i := (u_i,v_i)$ is 
an edge of $G - {\cal E}(P)$.
The vertices $u_i,\ i \in {\bf Z}_5$ are distinct 
(a trivial consequence
of the fact that $G$ has no 3-cycles or 4-cycles).
Suppose $G - {\cal E}(P)$ can be edge-3-colored, and 
$\gamma \in EC(G - {\cal E}(P))$. \medskip  

   (A) By the ``Parity Lemma'' (the special case of it
in Remark 2.12(d)), one has that \break 
(i) for exactly three distinct indices 
$q,r,s \in {\bf Z}_5$, the equality 
$\gamma(\epsilon_q) = \gamma(\epsilon_r) = 
\gamma(\epsilon_s)$ holds, and
(ii) for the other two indices $t$ and $u$ in 
${\bf Z}_5$, the colors
$\gamma(\epsilon_t)$, $\gamma(\epsilon_u)$, and 
$\gamma(\epsilon_q)$ 
($=\gamma(\epsilon_r)=\gamma(\epsilon_s)$) are 
distinct. \medskip

   (B) In (A), if $\{q,r,s\} = \{k-1,k,k+1\}$ for some $k \in {\bf Z}_5$
(that is, if $q,r,s$ are in some order ``consecutive'' in ${\bf Z}_5$),
then $\gamma$ extends uniquely to an edge-3-coloring of the entire
cubic graph $G$ itself. \medskip

   (C) In (A), if instead $\{q,r,s\} = \{k-2,k,k+2\}$ for 
some $k \in {\bf Z}_5$ (recall Remark 4.2(A)),
then $\gamma$ does not extend to an edge-3-coloring of $G$. \medskip

   (D) Paragraph (A) above has the following equivalent formulation:
For any given $\delta \in ED(G - {\cal E}(P))$, three of 
the edges
$\epsilon_i$, $i \in {\bf Z_5}$ belong to the same class 
(of the three classes in the decomposition $\delta$), 
and the other two of the edges $\epsilon_i$
belong respectively to the other two classes. \hfil\break

{\bf Remark 4.4.}\ \ Theorem 4.5 below will be motivated 
here by the following ``historical'' information: \medskip

   (A) Consider a simple cubic {\it planar\/} graph $H$ 
which has no ``bridge'' (an edge whose removal would
split the remainder of $H$ into two disjoint pieces). 
The Four Color Theorem states that there exists a 
``face-4-coloring'' of $H$ --- a coloring $\Gamma$ of the 
{\it faces\/} of $H$ with (say) the four colors 
$0, a, b, c$ (the elements of ${\bf Z}_2 \times {\bf Z}_2$ 
--- see (2.1)) such that no two
contiguous faces are assigned the same color. 
It is well known that any face-4-coloring $\Gamma$ of $H$ 
with the colors $0, a, b,c$ induces an edge-3-coloring 
$\gamma$ of $H$ with the colors $a,b,c$ in which for each 
edge $e$ of $H$,
$$ \gamma(e) = \Gamma(F_1) + \Gamma(F_2)     \eqno (4.1) $$
where $F_1$ and $F_2$ are the faces of $H$ that share the 
edge $e$ as a common border.  
(Recall that if $x$ and $y$ are any two distinct
elements of ${\bf Z}_2 \times {\bf Z}_2$ --- either one 
can be 0 --- then $x+y \neq 0$.)\ \ 
It is also well known that conversely, for any 
$\gamma \in EC(H)$, any face $F$ of $H$, and any element
$x \in {\bf Z}_2 \times {\bf Z}_2$, there exists a (unique) 
face-4-coloring $\Gamma$ of $H$ such that $\Gamma(F) = x$ 
and (4.1) holds for every edge $e$ of $H$.
Such connections between face-4-colorings and
edge-3-colorings go back to the work of Peter Tait 
in 1880 (alluded to after Definition 2.4),
about a century before the Four Color Theorem was proved; 
see [Wi, Chapter 6].   
\medskip

   (B) In a paper published in 1879, Alfred Kempe  
gave an intended  ``proof'' of the (then not yet proved) 
Four Color Theorem.  
The key facet of his argument was his claim that if $H$ 
is a simple cubic planar graph with girth 5 and no 
``bridge'', $P$ is a pentagon in $H$, and $\Gamma$ is a face-4-coloring of all faces of $H$ except the face $P$, 
then (if necessary) after successive interchanges of colors 
along two different connected two-color regions (now known 
as ``Kempe chains''), the coloring $\Gamma$ can be 
extended to include the face $P$ and thereby produce a 
face-4-coloring of (all of) $H$.
Eleven years later, in a paper published in 1890, 
Percy Heawood showed with an example that Kempe's intended ``proof'' --- in particular, the key facet of it described above --- is invalid.
(Heawood's example was of course not a counterexample to 
the Four Color Theorem itself; it just exposed a fatal 
flaw in Kempe's intended argument for it.)\ \ 
For a detailed description of Kempe's wrong ``proof'' and 
of Heawood's example, see [Wi, Chapters 5--7].
What (intentionally, by design) ``goes wrong'' in Heawood's example, will be referred to here informally as ``Heawood's monkey wrench''.
\medskip

   (C) Adapting the connection between face-4-colorings 
and edge-3-colorings (for simple cubic planar graphs
with no ``bridge'') 
described in paragraph (A) above, one can transcribe the description of ``Heawood's monkey wrench'' from its 
original form (involving face-4-colorings)
into a form involving edge-3-colorings.
Such an ``edge-3-coloring'' version of
``Heawood's monkey wrench'' can then be adapted
to some simple cubic graphs that are not planar.
\medskip

   (D) In particular, suppose $G$ is a snark.  
Of course, by the Four Color Theorem itself (and 
paragraph (A) above), $G$ is nonplanar.
If $P$ is a pentagon in $G$, then for the 
quasi-cubic graph $G-{\cal E}(P)$, if it is colorable, 
the ``edge-3-coloring'' version of 
``Heawood's monkey wrench'' inevitably 
{\it has\/} to occur, over and over and over again,
simply as a consequence of the fact that $G$ itself cannot 
be edge-3-colored.
The next theorem is based on that fact, and compiles some information that follows naturally from it. \hfil\break

   {\bf Theorem 4.5.}\ \ {\sl Suppose $G$ is a snark.

   (A) If $H$ is a connected union of pentagons in $G$ 
(or $H$ is simply a single pentagon in $G$), then the 
numbers $\psi(G,e),\ e \in {\cal E}(H)$ are equal.

   (B) Suppose $P$ is a pentagon in $G$ with vertices 
$v_i$ and edges $(v_i, v_{i+1})$, $i \in {\bf Z}_5$ 
(recall Convention 4.1).  
For each $i \in {\bf Z}_5$, let $\epsilon_i$ denote the 
edge of $G - {\cal E}(P)$ that is connected to the vertex 
$v_i$.
Refer to conclusion (A) above, and also to Remark 4.3(D).
Then the following two statements (i), (ii) hold:

(i) ${\rm card}\thinspace ED(G-{\cal E}(P)) = 5 \cdot \psi(G,e)$ 
where $e$ is any edge of $P$.

(ii) For each $i \in {\bf Z}_5$, 
${\rm card}\thinspace \{\delta \in ED(G-{\cal E}(P)):
\epsilon_{i-2} \sim \epsilon_i \sim \epsilon_{i+2} \}\ = \ 
\psi(G,e)$ where $e$ is any edge of $P$.} \hfil\break

   In conclusion (A), it is understood that the union of 
{\it all\/} pentagons in $G$ may be a disconnected subgraph 
of $G$, and that the common number $\psi(G,e)$ for the edges $e$ of one connected component of that subgraph may be different from the common number $\psi(G,e)$ for the
edges $e$ of another connected component.

   Of course in conclusion (B)(i)(ii), in the phrase 
``where $e$ is any edge of $P$'', $e$ is of course 
$(v_i, v_{i+1})$ for some (any) $i \in {\bf Z}_5$. 
\medskip

   Theorem 4.5, and its proof given below, are taken
from [Br1, Theorem 2 and Corollary 1], but with
some changes in terminology and style.
  
   It was also pointed out there (in [Br1, Theorem 2]) 
that in the context of Theorem 4.5(B) here,
as a simple consequence of conclusion (ii), for any two
non-adjacent edges $d$ and $f$ of $P$, one has that 
${\rm card}\thinspace ED((G_d)_f) = \psi(G,e)$ for (any) 
$e \in {\cal E}(P)$.
\medskip 

   As manifested in the proof given here,
Theorem 4.5 is a by-product of repeated 
applications of ``Heawood's monkey wrench'' 
in its edge-3-coloring form 
(see Remark 4.4(B)(C)(D) again).
In the argument below, that form of
``Heawood's monkey wrench'' plays its main role
in Lemma 3 and its proof.
\hfil\break

   In a slightly different way, Theorem 4.5 can also 
be derived, with a little work, as a corollary or 
by-product of some arguments of 
K\'aszonyi [K\'a3, Section 3].
(Note that if $G$ is a snark with a pentagon $P$,
and $e$ is an edge of $P$, then --- recall
Theorem 3.3(B) --- the two orthogonal edges 
$d_1$ and $d_2$ of $G_e$ in the context of 
Notations 3.2 are ``near'' to each other in the 
terminology of [K\'a3, Section 3].) 
\hfil\break

   {\bf Proof.}\ \ The proof of Theorem 4.5 will first 
proceed through a series of definitions, lemmas, etc.\  (numbered in order as 0, 1, 2, \dots, 9)
that will establish statement (B), together with 
statement (A) for the special case of a single pentagon.
After that, at the very end, a final argument will be given 
to establish statement (A) in its full generality.

   Throughout this proof, $G$ is a given snark.  
It is assumed to have at least one pentagon, say $P$,
(consistent with the hypotheses of each of statements
(A) and (B)).
However, the quasi-cubic graph $G - {\cal E}(P)$
is not assumed to be colorable; definitions and lemmas
below are allowed to be ``vacuous''. \medskip

   {\bf Proof of statement (B) (and of statement (A)
for one pentagon).}\ \ The proof will start 
with a ``Context'' that will provide the setting for the 
entire argument for statement (B) (and for statement (A)
for one pentagon). \medskip

   {\bf Context 0.}\ \ Suppose $P$ is a pentagon
(5-cycle) in $G$.
Let the vertices of $P$ be denoted in order as
$v_0, v_1, v_2, v_3, v_4$, with the edges of $P$ being
$(v_i, v_{i+1})$  
(of course using Convention 4.1). 
For each $i \in \{0,1,2,3,4\}$, let $\epsilon_i$ denote 
the edge of $G - {\cal E}(P)$ that is connected to the
vertex $v_i$, and let $w_i$ denote the ``other'' 
endpoint vertex of $\epsilon _i$, so that
$\epsilon_i = (w_i,v_i)$.
As a consequence of the snark $G$ having 
no 3-cycles or 4-cycles
(see Definition 2.9), the vertices $w_i$,
$i \in \{0,1,2,3,4\}$ (as well as the vertices
$v_i$, $i \in \{0,1,2,3,4\}$) are distinct. 
\medskip

   {\bf Definition 1.}\ \ Refer to Remark 4.3(A)(B)(C).
For each $k \in \{0,1,2,3,4\}$ (recall Convention 4.1),
let $Q_k$ denote the set of all 
$\gamma \in EC(G - {\cal E}(P))$ such that
$\gamma(\epsilon_{k-2}) = \gamma(\epsilon_k)
= \gamma(\epsilon_{k+2})$.

   It is understood (see Remark 4.3(A)) that  
for each $k \in \{0,1,2,3,4\}$ and each
$\gamma \in Q_k$, the colors
$\gamma(\epsilon_{k-1})$, $\gamma(\epsilon_k)$, and
$\gamma(\epsilon_{k+1})$ are distinct.

   By Remark 4.3(A)(B)(C) and the assumption that $G$ 
is a snark,
$$ EC(G - {\cal E}(P))\ =\  
\bigcup_{k \in \{0,1,2,3,4\}} Q_k   \eqno (4.2) $$
and (by a trivial argument) the sets $Q_k$ are (pairwise) disjoint. 
\medskip

   {\bf Definition 2.}\ \ In Context 0, define the mappings
$M : EC(G - {\cal E}(P)) \to EC(G - {\cal E}(P))$
and $M^* : EC(G - {\cal E}(P)) \to EC(G - {\cal E}(P))$
as follows.
The two definitions will be given together in three 
``steps'':
\medskip

   Suppose $\beta \in EC(G - {\cal E}(P))$.
   
   (i) Referring to (4.2) and the phrase right after it,
let $\kappa \in \{0,1,2,3,4\}$ be the index such that
$\beta \in Q_\kappa$.

   (ii) Referring to Definition 1 and Remark 4.3(A),
let the three distinct colors $s_-$, $s_o$, and $s_+$
(some permutation of the colors $a$, $b$, and $c$ in (2.1))
be defined by $s_o := \beta(\epsilon_{\kappa-2}) = \beta(\epsilon_\kappa) = \beta(\epsilon_{\kappa+2})$ and 
$s_- := \beta(\epsilon_{\kappa-1})$ and
$s_+ := \beta(\epsilon_{\kappa+1}).$

   (iii)(a) Let $M\beta$ denote the element of
$EC(G - {\cal E}(P))$ that one obtains from $\beta$
by interchanging the colors $s_o$ and $s_-$ along the
$s_os_-$-Kempe chain containing the edge
$\epsilon_{\kappa-1}$.

   (iii)(b) Let $M^*\beta$ denote the element of 
$EC(G - {\cal E}(P))$ that one obtains from $\beta$
by interchanging the colors $s_o$ and $s_+$ along the
$s_os_+$-Kempe chain containing the edge
$\epsilon_{\kappa+1}$.
\medskip

   {\bf Lemma 3.}\ \ {\sl In Context 0, suppose 
$\gamma \in EC(G - {\cal E}(P))$.
Referring to (4.2) and the phrase after it,
let $k \in \{0,1,2,3,4\}$ be the index such that
$\gamma \in Q_k$.
Then (recall Convention 4.1) 
the following four statements hold: \hfil\break
(i) $M\gamma \in Q_{k+2}$; \hfil\break
(ii) $M^*\gamma \in Q_{k-2}$; \hfil\break
(iii) $M^*M\gamma = \gamma$; and \hfil\break 
(iv) $MM^*\gamma = \gamma$.}
\medskip

   {\bf Proof.}\ \ The proofs of (ii) and (iv)
are respectively exactly analogous to
(are ``mirror images'' of) the proofs of (i) and (iii).
It will suffice to give the argument for (i) and (iii).
Here the arguments for (i) and (iii) will be given 
together.
The remaining paragraphs in this proof will be
labeled (P1), (P2), etc.

   (P1) As in the hypothesis, suppose
$\gamma \in EC(G - {\cal E}(P))$.
Let $k \in \{0,1,2,3,4\}$ be the index such that
$\gamma \in Q_k$. 
Let the (distinct) colors $x$, $y$, and $z$ be
defined by 
$$ x := \gamma(\epsilon_{k-2}) = \gamma(\epsilon_k)
= \gamma(\epsilon_{k+2}) \quad {\rm and} \quad
y := \gamma(\epsilon_{k-1}) \quad {\rm and} \quad
z := \gamma(\epsilon_{k+1}).   \eqno (4.3) $$

   (P2) Referring to (4.3), let $K$ be the $xy$-Kempe 
chain (for the coloring $\gamma$) containing the edge
$\epsilon_{k-1}$.
By Definition 2, the edge-3-coloring $M\gamma$ of
$G - {\cal E}(P)$ is obtained from the coloring $\gamma$
by the interchanging of the colors $x$ and $y$ along
that Kempe chain $K$.

   (P3) Recall that in the graph $G - {\cal E}(P)$, the 
vertex $v_{k-1}$ is univalent and is connected to the 
edge $\epsilon_{k-1}$.
It follows (recall Remark (d) in Definition 2.11)
that $K$ is a path (i.e.\ not a cycle) with two
endpoints, one of which is $v_{k-1}$ and the other is
some other univalent vertex. 
Since the only univalent vertices in $G - {\cal E}(P)$
are the five vertices $v_i$, $i \in \{0,1,2,3,4\}$,
the other endpoint of $K$ must be one of those 
vertices $v_i$ other than $v_{k-1}$.
That other endpoint cannot be $v_{k+1}$, because the
edge connected to it, namely $\epsilon_{k+1}$, is
colored $z$ (not $x$ or $y$) by $\gamma$ --- see (4.3).
Hence that other endpoint must be either
$v_k$, $v_{k-2}$, or $v_{k+2}$.

   (P4) If the other endpoint of $K$ (besides $v_{k-1}$) 
were $v_{k+2}$, then the mapping $M\gamma$ would assign 
the color $x$ to the edges $\epsilon_k$, 
$\epsilon_{k-1}$, and $\epsilon_{k-2}$, the color
$y$ to $\epsilon_{k+2}$, and the color $z$ to
$\epsilon_{k+1}$, and by Remark 4.3(B), the mapping
$M\gamma$ would extend to an edge-3-coloring of the
original graph $G$, contradicting the assumption that
$G$ is a snark.

   (P5) If the other endpoint of $K$ (besides $v_{k-1}$) 
were $v_k$, then the mapping $M\gamma$ would assign the
color $x$ to the edges $\epsilon_{k-1}$, 
$\epsilon_{k-2}$, and 
$\epsilon_{k+2} (= \epsilon_{k-3})$, the color
$y$ to $\epsilon_k$, and the color $z$ to
$\epsilon_{k+1}$, and in this case too by 
Remark 4.3(B), the mapping
$M\gamma$ would extend to an edge-3-coloring of the
original graph $G$, contradicting the assumption that
$G$ is a snark.

   (P6) Consequently, the other endpoint of $K$ (besides
$v_{k-1}$) has to be $v_{k-2}$.

   (P7) The coloring $M\gamma$ assigns the following
colors to the edges $\epsilon_i$, $i \in \{0,1,2,3,4\}$:
\break
$(M\gamma)(\epsilon_{k-1}) = (M\gamma)(\epsilon_k)
= (M\gamma)(\epsilon_{k+2}) = x$,
$(M\gamma)(\epsilon_{k-2}) = y$, and 
$(M\gamma)(\epsilon_{k+1}) = z$.
Using Convention 4.1, let us display this information
in a slightly different way:
$$
\eqalignno{ 
x &= (M\gamma)(\epsilon_k) = (M\gamma)(\epsilon_{k+2})
= (M\gamma)(\epsilon_{k+4}) \quad {\rm and} \cr
z &= (M\gamma)(\epsilon_{k+1}) \quad {\rm and} \quad
y = (M\gamma)(\epsilon_{k+3}).  & (4.4) \cr
}$$ 
By (4.4) and Definition 1, $M\gamma \in Q_{k+2}$.
Thus statement (i) in Lemma 3 holds.

   (P8) Our remaining task is to prove statement (iii) 
in Lemma 3. 
By (4.4) and Definition 2, since $M\gamma \in Q_{k+2}$,
the edge-3-coloring $M^*M\gamma$ is obtained from 
$M\gamma$ by the interchanging of the colors $x$ and $y$
along the $xy$-Kempe chain $K_1$ 
(for the coloring $M\gamma$) containing the edge 
$\epsilon_{(k+2)+1} = \epsilon_{k+3}
= \epsilon_{k-2}$.
(Here the coloring $\beta$, the index $\kappa$, and the 
colors $s_o$ and $s_+$ in Definition 2 are 
the coloring $M\gamma$, the index $k+2$, and the colors 
$x$ and $y$ in (4.4).)

   (P9) For comparison to the Kempe chain $K_1$
(for $M\gamma$) in paragraph (P8), note that the 
$xy$-Kempe chain $K$ (for the original coloring 
$\gamma$) in paragraph (P2) has the following two 
properties:   
(i) $K$ contains the edge edge $\epsilon_{k-2}$; and 
(ii) $K$ is an $xy$-Kempe chain for the coloring $M\gamma$.
Property (i) holds by paragraph (P6); and property (ii)
holds by paragraph (P2) and Remark (e)(i) in 
Definition 2.11,
By paragraph (P8), $K_1$ also satisfies (i) and (ii).
Hence $K_1$ is identical to $K$, by 
the first sentence after Definition 2.11(b). 
Thus by paragraph (P8), the coloring
$M^*M\gamma$ is obtained from $M\gamma$ by the
interchanging of the colors $x$ and $y$ along $K$.
Hence by paragraph (P2) and Remark (e)(ii) in 
Definition 2.11, statement (iii) in Lemma 3 holds.
This completes the proof of Lemma 3.
\medskip

   {\bf Lemma 4.}\ \ {\sl Refer to Definition 1.
In Context 0, for each $k \in \{0,1,2,3,4\}$,
${\rm card}\thinspace Q_k = 
{\rm card}\thinspace Q_{k+2}$.}
\medskip

   {\bf Proof.}\ \ Suppose $k \in \{0,1,2,3,4\}$.
Refer to Definition 2.
Let us restrict the mapping $M$ in Definition 2 to
the domain $Q_k$, and let us also restrict the
mapping $M^*$ in Definition 2 to the domain
$Q_{k+2}$.
Under these restrictions, by Lemma 3(i)(ii), we thereby
have that $M$ maps $Q_k$ into $Q_{k+2}$, and
$M^*$ maps $Q_{k+2}$ into $Q_k$.
By Lemma 3(iii) and the usual trivial argument, 
the (restricted) mapping $M$ is one-to-one (as a
mapping of $Q_k$ into $Q_{k+2}$).
By Lemma 3(iv) and the usual trivial argument,
the (restricted) mapping $M$ is also onto (as a
mapping of $Q_k$ into $Q_{k+2}$).
Hence the (restricted) mapping $M$ gives a one-to-one
correspondence between the sets $Q_k$ and $Q_{k+2}$.
Lemma 4 follows.
\medskip

   {\bf Step 5.}\ \ In Context 0, define the nonnegative integer $J := {\rm card}\thinspace Q_0$.
Then by four applications of Lemma 4,
$$ J\ =\ {\rm card}\thinspace Q_0\ =\ 
{\rm card}\thinspace Q_2\ =\ {\rm card}\thinspace Q_4\ =\ 
{\rm card}\thinspace Q_1\ =\ {\rm card}\thinspace Q_3.
\eqno (4.5) $$
Hence by eq.\ (4.2) and the phrase right after it,
$$ {\rm card}\thinspace EC(G - {\cal E}(P))\ =\ 5J.
\eqno (4.6) $$  

   {\bf Step 6.}\ \ In Context 0, suppose $e$ is any
edge of $P$, and let $k \in \{0,1,2,3,4\}$ denote the 
index such that $e = (v_{k-2}, v_{k+2})$ (recall
Convention 4.1).

   If $\mu$ is any edge-3-coloring of the cubic
graph $G_e$, then trivially there is a (unique)
edge-3-coloring $\gamma$ of 
$G - {\cal E}(P)$ that meets the following conditions
(referring to vertices $w_i$ in Context 0):
$$ \eqalignno{ 
\gamma(\epsilon_{k-2}) &= \mu((w_{k-2}, v_{k-1})), \quad \gamma(\epsilon_{k+2}) = \mu((w_{k+2}, v_{k+1})), \quad
{\rm and} \cr
 \gamma(\varepsilon) &= \mu(\varepsilon)\ 
{\rm for\ every\ edge}\ \varepsilon\ {\rm of}\ 
G-{\cal E}(P)\ {\rm other\ than}\ 
\epsilon_{k-2}\ {\rm and}\ \epsilon_{k+2}.  & (4.7) \cr
} $$

   {\bf Lemma 7.}\ \ {\sl In Context 0, suppose $e$ is any
edge of $P$, and let $k \in \{0,1,2,3,4\}$ denote the 
index such that $e = (v_{k-2}, v_{k+2})$ (recall
Convention 4.1).

   (i) For any $\gamma \in Q_{k-2} \cup Q_k \cup Q_{k+2}$,
there exists a unique edge-3-coloring $\mu$ of the
cubic graph $G_e$ such that (4.7) holds.

   (ii) For any $\gamma \in Q_{k-1} \cup Q_{k+1}$,
there does not exist an edge-3-coloring $\mu$ of 
$G_e$ such that (4.7) holds.}
\medskip 

   {\bf Proof.}\ \ Let us prove (ii) first.
Consider first the case where $\gamma \in Q_{k-1}$.  
Suppose there were to exist an edge-3-coloring 
$\mu$ of $G_e$ such that (4.7) holds.
Then one would have 
$ \gamma(\epsilon_{k-1}) = \gamma(\epsilon_{k+1}) =
\gamma(\epsilon_{k+2})$.
To the second and third terms there, one can apply (4.7), 
and one obtains 
(again recall the vertices $w_i$ in Context 0) 
$$\mu(\epsilon_{k+1}) = \gamma(\epsilon_{k+1}) =
\gamma(\epsilon_{k+2}) = \mu((w_{k+2}, v_{k+1})).
\eqno (4.8) $$
However, since the vertex $v_{k+1}$ is an endpoint of 
both of the edges $(w_{k+2}, v_{k+1})$ and 
$\epsilon_{k+1}$ (in $G_e$), one must have 
$\mu(\epsilon_{k+1}) \neq \mu((w_{k+2}, v_{k+1}))$, 
which contradicts (4.8).
Thus (if $\gamma \in Q_{k-1}$) there cannot exist an
edge-3-coloring $\mu$ of $G_e$ such that (4.7) holds.

   By an analogous (``mirror image'') argument, 
one has that if $\gamma \in Q_{k+1}$, there cannot exist
an edge-3-coloring $\mu$ of $G_e$ such that (4.7) holds.
That completes the proof of statement (ii) in Lemma 4.7.
\medskip

   {\it Proof of statement (i).}\ \ For a given 
edge-3-coloring $\gamma$ of $G - {\cal E}(P)$
such that $\gamma \in Q_{k-2} \cup Q_k \cup Q_{k+2}$,
in order to ``extend'' it to an edge-3-coloring of
$G_e$ --- more precisely, in order to define an edge-3-coloring $\mu$ of $G_e$ that satisfies (4.7) --- one would need to assign colors $\mu(.)$ to the two remaining edges 
$(v_{k-1}, v_k)$ and $(v_k, v_{k+1})$ (of $G_e$), 
in such as way as to avoid giving 
the same color to two adjacent edges.

   Consider first the case where $\gamma \in Q_k$.
Let the (three distinct) colors $x$, $y$, and $z$ be
as in (4.3) (in the proof of Lemma 3).  
There exists a unique edge-3-coloring of $G_e$ such
that (4.7) holds.
It is obtained by assigning the following
colors to the two remaining edges:
$\mu((v_{k-1},v_k)) := z$ and
$\mu((v_k,v_{k+1})) := y$.
 
   Next consider the case where $\gamma \in Q_{k-2}$.
Let the (three distinct) colors $x$, $y$, and $z$ be
defined by
$$ x := \gamma(\epsilon_{k+1}) = \gamma(\epsilon_{k-2})
= \gamma(\epsilon_k) \quad {\rm and} \quad
y := \gamma(\epsilon_{k+2}) \quad {\rm and} \quad
z := \gamma(\epsilon_{k-1}). $$
This is simply a version of (4.3) with $-2$ added to each 
index. In this case, there is a unique edge-3-coloring
$\mu$ of $G_e$ that satisfies (4.7).  
It is obtained by assigning the following colors
to the two remaining edges:
$\mu((v_{k-1},v_k)) := y$ and
$\mu((v_k,v_{k+1})) := z$.

   For the remaining case $\gamma \in Q_{k+2}$ the
argument is exactly analogous to (is a ``mirror image'' 
of) the argument for the case $\gamma \in Q_{k-2}$.
That completes the proof of statement (i), and of Lemma 7.
\medskip         

   {\bf Lemma 8.}\ \ {\sl In Context 0, for every 
edge $e$ of the pentagon $P$, one has that}
$$ {\rm card}\thinspace EC(G_e)\ =\ 3J.  \eqno (4.9) $$
   \indent{\bf Proof.}\ \ By Step 6 and Lemma 7, eq.\ (4.7)
gives a one-to-one correspondence between the set of all
edge-3-colorings $\mu$ of $G_e$ and the set of all
$\gamma \in Q_{k-2} \cup Q_k \cup Q_{k+2}$.
Hence by (4.5) (and the phrase right after (4.2)),
eq.\ (4.9) holds.
\medskip

   {\bf Step 9.}\ \ In this step, it will be
convenient to slightly abbreviate the earlier notation 
$\{\delta \in ED(G - {\cal E}(P)): 
\epsilon_{k-2} \sim \epsilon_k \sim \epsilon_{k+2}\}$
to simply 
$\{ED(G - {\cal E}(P)): 
\epsilon_{k-2} \sim \epsilon_k \sim \epsilon_{k+2}\}$.
   
   Refer again to Context 0.
By eq.\ (2.4) in Remark 2.6, one has that the integer 
${\rm card}\thinspace EC(G-{\cal E}(P))$ is
a multiple of 6.
It now follows from (4.6) that the (nonnegative) integer 
$J$ must be a multiple of 6.
Define the nonnegative integer $L$ by $L := J/6$.
Then by (4.6) and (4.9), together with (again)
eq.\ (2.4), 
$$ {\rm card}\thinspace ED(G - {\cal E}(P))\ =\ 5L;
\eqno (4.10) $$
and 
$${\rm card}\thinspace ED(G_e) = 3L \quad
{\rm for\ each\ edge}\ e\ {\rm of\ the\ pentagon}\ P.   
\eqno (4.11) $$

   Also, for any given $k \in \{0,1,2,3,4\}$, 
any $\delta \in 
\{ED(G - {\cal E}(P)): 
\epsilon_{k-2} \sim \epsilon_k \sim \epsilon_{k+2}\}$ 
gives rise to exactly 6 edge-3-colorings $\gamma \in Q_k$
(with the edges in any of the three classes in the
decomposition $\delta$ being given the same color),
since there are exactly 6 permutations of the three
colors $a$, $b$, and $c$.
Of course any such coloring $\gamma$ arises from exactly 
one such decomposition $\delta$.
Hence for each $k \in \{0,1,2,3,4\}$,
${\rm card}\thinspace Q_k = 6 \cdot 
{\rm card}\thinspace \{ED(G - {\cal E}(P)): 
\epsilon_{k-2} \sim \epsilon_k \sim \epsilon_{k+2}\}$. 
Hence by (4.5),
$$ \forall k \in \{0,1,2,3,4\}, \quad
{\rm card}\thinspace \{ED(G - {\cal E}(P)): 
\epsilon_{k-2} \sim \epsilon_k \sim \epsilon_{k+2}\} = L.
\eqno (4.12) $$

   Now by (3.7) and (4.11),
$$ \psi(G,e) = L \quad 
{\rm for\ each\ edge}\ e\ {\rm of\ the\ pentagon}\ P.   
\eqno (4.13) $$
One now obtains sub-statement (i) in statement (B)
(in Theorem 4.5) by substituting (4.13) into (4.10),
and one obtains sub-statement (ii) in statement (B)
by substituting (4.13) into (4.12).
That, together with (4.13) itself, completes the proof of statement (B) (and of statement (A) for one pentagon) in Theorem 4.5.
\medskip

   {\bf Proof of statement (A).}\ \ Recall that in the
case where $H = P$ itself for some pentagon $P$ in $G$,
from (4.13) in the proof above (for statement (B)
and for statement (A) for this particular pentagon $P$),
one already has established
that the numbers $\psi(G,e)$, $e \in P$ are equal.
This special case will be tacitly used below.

   Now suppose instead that $H$ is a connected union of 
two or more pentagons in $G$.
Suppose $e_1$ and $e_2$ are any two distinct edges of $H$.
It suffices to prove that
$$ \psi(G,e_1) = \psi(G, e_2).  \eqno (4.14) $$

   Since $H$ is connected, there is a finite sequence of 
edges $\varepsilon_0, \varepsilon_1, \dots, \varepsilon_n$
in $H$ such that $\varepsilon_0 = e_1$,
$\varepsilon_n = e_2$, and 
for each $i \in \{0,1,\dots, n-1\}$, the edges 
$\varepsilon_i$ and $\varepsilon_{i+1}$ are adjacent. 
If one can show that 
$\psi(G, \varepsilon_i) = \psi(G, \varepsilon_{i+1})$
for each $i \in \{0,1,\dots, n-1\}$, then (4.14) will
follow.
Hence, it suffices to prove (4.14) for the case
where $e_1$ and $e_2$ themselves are adjacent.

   Suppose $e_1$ and $e_2$ are adjacent.
If they belong to the same pentagon in $H$, then
from the comments above (involving (4.13)), 
$\psi(G,e_1) = \psi(G,e_2)$ 
and we are done.
Therefore, suppose instead that $e_1$ and $e_2$ do
{\it not\/} belong to the same pentagon in $H$.

   Let $v$ denote the common endpoint vertex of the two
edges $e_1$ and $e_2$.
Let $e_3$ be the third edge connected to $v$.
Now (by hypothesis) $e_1$ belongs to some pentagon $P_1$
in $H$.
This pentagon $P_1$ does not contain the edge $e_2$.
It follows that $P_1$ is forced to contain the edge $e_3$
(since every vertex in $P_1$, including $v$, is connected to
two edges in $P_1$).
Hence $\psi(G,e_1) = \psi(G,e_3)$.
Similarly, $e_2$ belongs to some pentagon $P_2$ in $H$,
$P_2$ does not contain $e_1$, hence $P_2$ must contain
$e_3$, and hence $\psi(G,e_2) = \psi(G,e_3)$.
Eq.\ (4.14) now follows.
That completes the proof of statement (A), and of
Theorem 4.5.
\hfil ////\break

{\bf Definition 4.6.}\ \ Refer to Theorem 4.5(A).  If $G$ 
is a snark and $H$ is a connected union of pentagons in $G$ (or $H$ is simply a single pentagon in $G$), then define the nonnegative integer
$$\psi(G,H) := \psi(G,e)       \eqno (4.15) $$ 
where $e$ is any edge of $H$.  \hfil\break

   The following well known procedure for combining two 
snarks, each having a pentagon, to create a ``bigger''
snark, was first used by 
Isaacs [Is, pp.\ 234-236] to create his ``double star'' 
snark from two disjoint copies of the flower snark 
snark $J_5$ (see Remark 2.14, and note that $J_5$
has a lone pentagon).       
This procedure also plays a role in the ``decomposition'' 
of a ``big'' snark with a (nontrivial) 5-edge cut set 
into two ``smaller'' snarks; 
see e.g.\ [CCW, Theorem 2]. \hfil\break

   {\bf Context 4.7.}\ \ Suppose $G'$ is a snark with a pentagon $P'$ with vertices $u_k$ and edges 
$(u_k,u_{k+1})$, $k \in {\bf Z}_5$ (recall Convention 4.1).
Let $G_0' := G' - \{u_k, k \in {\bf Z}_5\}$ (the graph
that one derives from $G'$ by deleting the five vertices
$u_k$ and all ten edges connected to them, including the
five edges in the pentagon $P'$).  
For each $k \in {\bf Z}_5$, let $t_k$ denote the vertex of $G_0'$ such that $(t_k,u_k)$ is an edge of $G'$.  
Since $G'$ is a snark (and hence has no 3-cycles or
4-cycles), these vertices $t_k,\ k \in {\bf Z}_5$ are distinct.

   Suppose $G^*$ is a snark that is disjoint from $G'$ and 
has a pentagon $P^*$ with vertices $v_k$ and edges 
$(v_k,v_{k+1})$, $k \in {\bf Z}_5$.
Let $G_0^* := G^* - \{v_k, k \in {\bf Z}_5\}$.  
For each $k \in {\bf Z}_5$, let $w_k$ denote the vertex of $G_0^*$ such that $(v_k,w_k)$ is an edge of $G^*$.  
These vertices $w_k,\ k \in {\bf Z}_5$ are distinct.

   Let $G$ denote the simple cubic graph that consists of $G_0'$, $G_0^*$, and the five new edges $(t_k, w_{2k})$, 
$k \in {\bf Z}_5$ (again recall Convention 4.1).
\hfil\break

   As is well known --- e.g.\ from Isaacs [Is] in 
connection with his ``double star'' snark alluded to above, 
right before Context 4.7 --- the graph $G$ in
Context 4.7 is a snark.  
Here let us quickly review the well known proof 
(from [Is, Theorem 4.2.1]) that the 
graph $G$ cannot be edge-3-colored.  
(To verify that $G$ is at least 4-edge-connected and
has girth at least 5, requires separate arguments, which
we shall not go into here.) \medskip

   Suppose (to seek a contradiction) $\gamma$ were an
edge-3-coloring of $G$.  
Then by Remark 2.12(d) 
(see Lemma 2.13 and adapt the Remark after it), 
one has that 
(i) for (exactly) three indices $q,r,s \in {\bf Z}_5$,
$\gamma((t_q, w_{2q})) = \gamma((t_r, w_{2r})) =
\gamma((t_s, w_{2s}))$, and
(ii) for the other two indices $i$ and $j$ in ${\bf Z}_5$,
the colors $\gamma((t_q, w_{2q}))$, $\gamma((t_i, w_{2i}))$,
and $\gamma((t_j, w_{2j}))$ are distinct.
If the three indices $q$, $r$, and $s$ are (in some order)
consecutive in ${\bf Z}_5$, then
$\gamma$ would induce an edge-3-coloring
of $G'$ (for each $k \in {\bf Z}_5$, the color
$\gamma((t_k, w_{2k}))$ is given to the edge
$(t_k, u_k)$ in $G'$, and then one applies
Remark 4.3(B)); but that would contradict the assumption 
that $G'$ is a snark.
If instead the three indices $q$, $r$, and $s$ are
not all consecutive in ${\bf Z}_5$ (in any order), then by Remark 4.2(C), the indices $2q$, $2r$, and $2s$ {\it are\/} (in some order) consecutive in ${\bf Z}_5$, and $\gamma$
would analogously induce an edge-3-coloring of $G^*$, contradicting the assumption that $G^*$ is a snark.
Thus an edge-3-coloring $\gamma$ of $G$ cannot exist.     
\hfil\break

   {\bf Theorem 4.8.}\ \ {\sl In Context 4.7, the following two statements hold (see Definition 4.6):\ \   
(A) For any edge $e$ of $G_0^*$, $e$ is an edge of
$G$ and
$$ \psi(G,e) = \psi(G^*,e) \cdot \psi(G',P')\ . 
\eqno (4.16) $$  
(B) For any edge $e$ of $G_0'$, $e$ is an edge of
$G$ and
$$ \psi(G,e) = \psi(G',e) \cdot \psi(G^*,P^*)\ . 
\eqno (4.17) $$}

   This theorem, and its proof given below, are due to the
author [Br2, Theorem 3.2 and the sentence after it].
The proof was given there somewhat tersely.  It will be 
repeated here in more generous detail. \hfil\break

   {\bf Proof.}\ \ It will suffice to give the argument for 
statement (A).
The proof for statement (B) is exactly analogous, and will
not be given explicitly here.
\medskip

   {\bf Proof of statement (A).}\ \ As in the hypothesis
of statement (A), suppose $e$ is an edge of $G_0^*$.
Then trivially from Context 4.7, $e$ is an edge of $G$.
For notational convenience, we shall carry out the proof
here for the case where $e$ is not connected to any of
the vertices $w_k$, $k \in {\bf Z}_5$.
(If instead $e$ were connected to one or two of those
vertices $w_k$, the proof would be essentially the same,
with only minor notational changes.)
\medskip

   Suppose $\gamma$ is a coloring of $G_e$.  
Then $\gamma$ induces the following coloring $\gamma'$
of the quasi-cubic graph $G'- {\cal E}(P')$:
$$ \eqalignno {\gamma'(d) := \gamma(d)\ \ 
&{\rm for}\ d \in G_0',\ \ {\rm and} \cr
\gamma'((t_k, u_k)) := \gamma((t_k, w_{2k}))\ \  
&{\rm for}\ k \in {\bf Z}_5. & (4.18) \cr   
}$$

   By Remark 2.12(d), applied to (say) $G'-{\cal E}(P')$, 
the colors $\gamma((t_k, w_{2k}))$, $k \in {\bf Z_5}$ 
are the same for three indices $k = q,r,s$ and distinct 
for the other two (i.e.\ with the other two colors each appearing on exactly one of those other two edges).
By Remark 4.3(A)(B) and the assumption that $G'$ is a 
snark, those three indices $q,r,s$ cannot all be 
(in any order) consecutive in ${\bf Z}_5$ (mod 5).
Hence by Remark 4.2(C), the indices $2q, 2r, 2s$ are
(in some order) all consecutive in ${\bf Z}_5$.
Hence by Remark 4.3(B), $\gamma$ induces a unique
coloring $\gamma^*$ of $G_e^*$ as follows:
$$ \eqalignno {\gamma^*(d) := \gamma(d)\ \ 
&{\rm for}\ d \in 
G_e^* - \{v_k, k \in {\bf Z_5}\},\ \  {\rm and} \cr
\gamma^*((v_{2k}, w_{2k})) := \gamma((t_k, w_{2k}))\ \ 
&{\rm for}\ k \in {\bf Z}_5. & (4.19) \cr   
}$$
(The colors $\gamma^*(d)$ for the edges $d$ of
$P^*$ are not specified here; they will be uniquely
determined by the colors $\gamma^*((v_{2k}, w_{2k}))$,
$k \in {\bf Z}_5$.)\ \ 
Thus the coloring $\gamma$ of $G_e$ induces an ordered
pair $(\gamma', \gamma^*)$ such that (see (4.18) and 
(4.19))
$$ \eqalignno{
&\gamma'\ {\rm is\ a\ coloring\ of}\ G' - {\cal E}(P), \cr
&\gamma^*\ {\rm is\ a\ coloring\ of}\ G_e^*,\ \  
{\rm and} \cr
&\gamma'((t_k, u_k)) = \gamma^*((v_{2k}, w_{2k}))\ \  
{\rm for}\ k \in {\bf Z}_5.  & (4.20) \cr 
} $$  

   Conversely, an ordered pair $(\gamma', \gamma^*)$ as
in (4.20) induces a (unique) coloring $\gamma$ of $G_e$
via (4.18) and (4.19).
Thereby one has a one-to-one correspondence between
colorings $\gamma$ of $G_e$ and ordered pairs
$(\gamma', \gamma^*)$ satisfying (4.20).

   Now suppose $\gamma^*$ is a coloring of $G_e^*$.
Then by Remark 4.3(A)(B)(C), there exists a permutation
$x,y,z$ of the colors ($a$, $b$, and $c$) and an element
$\ell$ of ${\bf Z}^5$ such that
$$
\gamma^*((v_k,w_k)) = \cases {
x & for $k \in \{\ell-1, \ell, \ell+1\}$ \cr
y & for $k = \ell-2$ \cr
z & for $k = \ell+2$. \cr}   \eqno (4.21)  
$$
Let $j \in {\bf Z}_5$ be defined by $j := 3\ell$ (mod 5)
(and hence $2j = 6\ell = \ell$ (mod 5)).
If $\gamma'$ is a coloring of the quasi-cubic graph
$G' - {\cal E}(P)$, then (the last line of) (4.20)
holds if and only if  
$$
\gamma'((t_k,u_k)) = \cases {
x & for $k \in \{j-2, j, j+2\}$ \cr
y & for $k = j-1$ \cr
z & for $k = j+1$. \cr}   \eqno (4.22)  
$$
By (4.15) and Theorem 4.5(B)(ii) (and Remark 4.3(D)), 
there are exactly
$\psi(G',P')$ 3-edge-decompositions of
$G'- {\cal E}(P')$ such that the edges
$(t_k,u_k)$, $k \in \{j-2,j,j+2\}$ are in the same
class (in the decomposition) and the edges
$(t_k,u_k)$, $k \in \{j-1, j+1\}$ are in the other
two classes respectively.
Hence (similarly to Remark 2.6)
there are exactly $\psi(G',P')$ colorings $\gamma'$ of
$G'- {\cal E}(P')$ such that (4.22) holds,
equivalently such that the ordered pair
$(\gamma', \gamma^*)$ is as in (4.20).
We have shown that this holds for an arbitrary
coloring $\gamma^*$ of $G^*_e$.

   Recall from Definition 3.4 that there are exactly
$18 \cdot \psi(G^*,e)$ colorings $\gamma^*$ of $G_e$.   
Hence by the preceding paragraph, there are exactly
$18 \cdot \psi(G^*,e) \cdot \psi(G',P')$ ordered
pairs $(\gamma', \gamma^*)$ satisfying (4.20), and hence
that is the number of colorings of $G_e$. 
Hence by Definition (3.4), eq.\ (4.16) holds.
That completes the proof of statement (A), and of
Theorem~4.8.
\hfil ////\break

\noindent {\bf 5.\ \ An example of $\psi(G,e)$ for a superposition} 
\hfil\break

   Isaacs [Is] presented some methods for creating 
arbitrarily large snarks.
One of those methods involved a particular procedure 
(a ``four-edge connection'', which Isaacs called
a ``dot product'') for ``combining'' two disjoint 
``smaller'' snarks to form a ``bigger'' snark.
(The formal definition of ``dot product'' can be found 
in [Is] or [Ga], and the procedure will be described just 
informally in comments after Context 5.2 below.
See also the comments in Section 7.2(B).)\ \ 
In general, in the ``dot product'', the roles of the two
``smaller'' snarks are not ``symmetric'' to each other.
However, in a certain class of special cases, the
roles of the two ``smaller'' snarks (when looked at in the right way) are in fact ``symmetric'' to each other. 
We shall allude to that type of special case 
here (as in [Br2]) as a ``symmetric dot product''.
\medskip

   If $G_1$ and $G_2$ are two disjoint snarks, and $G$
is a snark obtained from $G_1$ and $G_2$ via a ``symmetric
dot product'', then for every edge $e$ of $G$, there is a 
``natural'' choice of edges $e_1$ of $G_1$ and 
$e_2$ of $G_2$ such that
$\psi(G,e)$ is either equal to 
$2 \cdot \psi(G_1,e_1) \cdot \psi(G_2,e_2)$ or equal to   
$3 \cdot \psi(G_1,e_1) \cdot \psi(G_2,e_2)$. 
\medskip

   That ``fact'' was given a precise formulation in 
[Br2, Theorem 2.2 and subsequent sentence].
We shall not elaborate further on it here.  
That ``fact'' (the precise formulation of it) was 
applied in an induction argument 
(using the Petersen graph and Theorem 3.5 as the 
starting point as well as in the induction step) 
to prove the following theorem 
([Br2, Theorem 1.4]): \hfil\break  

   {\bf Theorem 5.1.}\ \ {\sl Suppose $j$ and $k$ are 
each a nonnegative integer.
Then there exists a snark $G$ and an edge $e$ of $G$
such that $\psi(G,e) = 2^j \cdot 3^k$.}
\hfil\break

   In that paper [Br2], the following open problem was 
implicitly posed:
{\it For precisely what positive integers $n$ do there 
exist a snark $G$ and an edge $e$ of $G$ such that 
$\psi(G,e) = n$?}
\medskip

   The purpose of this paper here is to promote research
on this problem and on other closely related
problems, such as those posed in Section 7. 
\medskip

   For the problem here, to make further progress beyond
Theorem 5.1, one will need to apply much broader and more
flexible techniques for creating ``big'' snarks from
``small'' ones than just the ``symmetric dot product''.
Such a broad, flexible technique was devised by
Kochol [Ko1], who referred to it under the name
``superpositions''.
That technique appears to have much promise for this
problem and related ones. 
\medskip

   Kochol's [Ko1] technique involves starting with a
snark, and ``replacing'' one or more of its edges
by snarks in certain ways
(together with ``replacing'' the end-point vertices of 
those edges by, say, quite general appropriate quasi-cubic 
graphs --- not necessarily snarks), to create a ``bigger'' 
snark.
For a particularly fascinating application of Kochol's 
technique, the reader is referred to the paper [Ko1] 
itself (starting with Theorem 1 there and its proof),  
where snarks with arbitrarily large girth are constructed. 
If the notations in that paper seem a little bewildering
at first, the diagrams in that paper tell the key features 
of that story very well --- in light of alert, thoroughly pervasive use of both Lemma 2.13 (the Parity Lemma) above 
and the properties of the Petersen graph.
\medskip 

Isaacs' ``dot product'', alluded to above, can be regarded
as the simplest example of a Kochol superposition.
(A little more on that below.)\ \ 
The rest of this section will be devoted to an illustration
of another simple Kochol superposition, in connection
with the numbers $\psi(G,e)$ for a given snark $G$ and
edge $e$ of $G$.
This particular superposition will provide a good simple 
illustration of, and motivation for, 
the use of the group arithmetic
on ${\bf Z}_2 \times {\bf Z}_2$ (as in Definition 2.1)
in the study of problems involving the 
numbers $\psi(G,e)$.
\hfil\break

   {\bf Context 5.2.}\ \ Suppose $G'$ is a snark, and
$E = (U,V)$ is an edge of $G'$.
Let $G_0' = G' - \{U,V\}$.
Suppose $T_{-2}, T_2$ are the vertices of $G_0'$ such that
$(T_i,U)$, $i \in \{-2,2\}$ are edges of $G'$.  
Suppose $W_{-2}, W_2$ are the vertices of $G_0'$ such that
$(V,W_i)$, $i \in \{-2,2\}$ are edges of $G'$.

   Suppose $G^*$ is a snark that is disjoint from $G'$.
Suppose $u$ and $v$ are (distinct) vertices of $G^*$ that
are not adjacent to each other in $G^*$.  
Suppose $u_i$, $i \in \{-1,0,1\}$ are the vertices of 
$G^*$ such that $(u, u_i)$ is an edge of $G^*$.
Suppose $v_i$, $i \in \{-1,0,1\}$ are the vertices of 
$G^*$ such that $(v_i, v)$ is an edge of $G^*$.
Let $G_0^* := G^* - \{u,v\}$.

   Let $G$ denote the simple cubic graph that consists of
$G_0'$, $G_0^*$, six new vertices $T_i$ and $W_i$,
$i \in \{-1,0,1\}$, and fourteen new edges
$(T_i, T_{i+1})$ and $(W_i, W_{i+1})$,
$i \in \{-2, -1, 0, 1\}$ and $(T_i, u_i)$ and
$(v_i, W_i)$, $i \in \{-1, 0, 1\}$.
\hfil\break

   (One of the vertices $u_i$ may be equal to one of 
the vertices $v_j$.
Any further equalities of those vertices is prohibited
by the requirement that the snark $G^*$ have
girth at least 5 --- recall Definition 2.9.)
\medskip

   Context 5.2 gives one way (certainly not the only one) 
to ``replace'' the edge $E$ in the snark $G'$ by the 
snark $G^*$ --- that is, within the snark $G'$, to ``superimpose'' the snark $G^*$ in place of the 
edge $E$.   
As part of the known information pertaining to Kochol's superpositions, it is well known that (under the
assumption that $G'$ and $G^*$ are each a snark) the
new graph $G$ is a snark.
As a review of the proof that $G$ cannot be colored,
one can simply carry out the proof of Theorem 5.3
below, but with $G_e$ and $G_e^*$ replaced by
$G$ and $G^*$; at (5.7) one would then obtain a coloring
of the snark $G^*$, a contradiction.
It takes extra arguments (which we shall not go into here)
to verify for $G$ the other properties in the definition
of ``snark'' (Definition 2.9). \medskip

   In the above superposition of the snark $G^*$ in place 
of the edge $E$, each of two (non-adjacent) vertices of 
$G^*$ were ``split'' into three ``strands'' that were
then ``hooked up to $G'$ where a vertex of $G'$ used to 
be''.
\medskip
  
   In other possible superpositions of $G^*$ in place 
of $E$, one might instead ``split'' each of two
nonadjacent {\it edges\/} of $G^*$ into two ``strands'' 
to be similarly reattached. 
In its simplest form, that is actually what is done in Isaacs' [Is] ``dot product'' alluded to 
above ---  producing a ``4-edge connection'' that 
``combines'' two snarks to form a bigger snark.
In a different but closely related context
(colorable simple cubic graphs with orthogonal edges), 
essentially the same procedure was also formulated
slightly earlier by K\'aszonyi [K\'a1]; 
see Section 7.2(B) in Section 7. 
\medskip   

   As yet another obvious alternative, one might 
``split'' one edge (of $G^*$) into two ``strands'' and one vertex (not connected to that edge) into three ``strands''.
\medskip

   In a superposition studied by McKinney [McK] (somewhat 
more complicated than that in Context 5.2), a pair of 
adjacent edges of a snark was ``replaced'' by a pair
of disjoint Petersen graphs.  
That will be discussed briefly in Section 6 below.
(That ``double superposition'' involved the ``splitting'' 
of both edges and vertices.)   
\hfil\break

   {\bf Theorem 5.3}\ \ {\sl In Context 5.2, suppose $e$ 
is an edge of $G_0^*$.
Then $e$ is an edge of $G$, and
$$ \psi(G,e) = 2 \cdot \psi(G^*,e) \cdot \psi (G',E).
\eqno (5.1) $$}

   {\bf Proof.}\ \ For convenience of notations, the
proof will be carried out in the case where the edge $e$ 
is not connected to any of the vertices $u_i$, $v_i$,
$i \in \{-1, 0,1\}$.  (In the case where $e$ is connected
to one of those vertices, the proof is essentially the
same, with only trivial changes in notations.)\ \ 
Refer to Definition 2.1.  
In the argument below, considerable use will be made of elementary properties of the group 
${\bf Z}_2 \times {\bf Z}_2$ (whose three
nonzero elements are the three colors).  

   Suppose $\gamma$ is a coloring of $G_e$. 
By Lemma 2.13 (and the Remark after it), 
$$ \gamma((T_{-2}, T_{-1})) + \gamma((T_1, T_2)) 
+ \sum_{k=-1}^1 \gamma((T_k, u_k)) = 0;  \eqno (5.2)$$
$$ \gamma((W_{-2}, W_{-1})) + \gamma((W_1, W_2)) 
+ \sum_{k=-1}^1 \gamma((v_k, W_k)) = 0;  \eqno (5.3)$$
and
$$ \sum_{k=-1}^1 \gamma((T_k, u_k))
+ \sum_{k=-1}^1 \gamma((v_k, W_k)) = 0. \eqno (5.4) $$
Since each element of ${\bf Z}_2 \times {\bf Z}_2$ is
its own inverse, it follows that 
for some $x \in {\bf Z}_2 \times {\bf Z}_2$,
$$ \eqalignno {
x\ &=\ \sum_{k=-1}^1 \gamma((T_k, u_k))\ 
=\ \sum_{k=-1}^1 \gamma((v_k, W_k)) \cr
&=\ \gamma((T_{-2}, T_{-1})) + \gamma((T_1, T_2))\ 
=\ \gamma((W_{-2}, W_{-1})) + \gamma((W_1, W_2)).
 & (5.5) \cr    
} $$
In (5.5), if $x \in \{a,b,c\}$ were to hold, then by 
(5.5) and Remark 2.2 (and the fact that in 
${\bf Z}_2 \times {\bf Z}_2$,
$x = y+z \Longleftrightarrow 0 = x + y + z$),
one would obtain a coloring
$\beta$ of the snark $G'$ (and hence a contradiction)
by defining $\beta(d) := \gamma(d)$ for 
$d \in {\cal E}(G_0')$, together with
$$ \eqalign{
\beta((T_{-2}, U)) &:= \gamma((T_{-2}, T_{-1})), \cr
\beta((T_2, U)) &:= \gamma((T_2, T_1)), \cr
\beta((W_{-2}, V)) &:= \gamma((W_{-2}, W_{-1})), \cr
\beta((W_2, V)) &:= \gamma((W_2, W_1)), \quad {\rm and} \cr
\beta(E) &:= x. \cr
} $$   
Hence $x = 0$ instead.
Thus by (5.5),
$$ \sum_{k=-1}^1 \gamma((T_k, u_k))\ 
=\ \sum_{k=-1}^1 \gamma((v_k, W_k))\ =\ 0.  \eqno (5.6) $$
Hence by Remark 2.2, $\gamma$ induces a coloring 
$\gamma^*$ of $G_e^*$ defined by
$$ \eqalignno{
\gamma^*(d) &:= \gamma(d) \quad {\rm for}\
d \in {\cal E}(G_e^* - \{u,v\}), \cr
\gamma^*((u,v_k)) &:= \gamma((T_k,u_k)) \quad {\rm for}\
k \in \{-1,0,1\}, \quad {\rm and} \cr
\gamma^*((v_k,v)) &:= \gamma((v_k,W_k)) \quad {\rm for}\
k \in \{-1,0,1\}. & (5.7) \cr 
} $$
Also, using (5.5) and the fact that $x = 0$, and then 
using a trivial extra coloring argument, one has that
$$ \eqalignno{
&\gamma((T_{-2}, T_{-1})) = \gamma((T_1, T_2)) 
= \gamma((T_0, u_0)), \quad {\rm and} \cr 
&\gamma((W_{-2}, W_{-1})) = \gamma((W_1, W_2))
= \gamma((v_0,W_0)). & (5.8) \cr
} $$  
Hence $\gamma$ induces a coloring $\gamma'$ of
$G_E'$ defined by   
$$ \eqalignno{
\gamma'(d) &:= \gamma(d) \quad {\rm for}\
d \in {\cal E}(G_0'), \cr
\gamma'((T_{-2},T_2)) &:= \gamma((T_0,u_0)), {\rm and} \cr
\gamma'((W_{-2},W_2)) &:= \gamma((v_0,W_0)). & (5.9) \cr 
} $$
Thus a given coloring $\gamma$ of $G_e$ induces an ordered
pair $(\gamma^*, \gamma')$ where (see (5.7) and (5.9))
$$ \eqalignno {
&\gamma^*\ {\rm is\ a\ coloring\ of}\ G_e^*, \cr
&\gamma'\ {\rm is\ a\ coloring\ of}\ G_E', \cr
&\gamma'((T_{-2}, T_2)) = \gamma^*((u,u_0)), \quad {\rm and} \cr
&\gamma'((W_{-2}, W_2)) = \gamma^*((v_0,v)). & (5.10) \cr
} $$
Conversely, if $(\gamma^*, \gamma')$ is an ordered pair satisfying (5.10), then it induces a unique coloring 
$\gamma$ of $G_e$ via (5.7) and (5.9).
(The colors $\gamma(d)$ for the remaining four edges
$d = (T_k, T_0), (W_k, W_0), 
\allowbreak k \in \{-1,1\}$ will trivially
be uniquely determined.)\ \ 
Thereby one obtains a one-to-one correspondence between
colorings $\gamma$ of $G_e$ and ordered pairs 
$(\gamma^*, \gamma')$ as in (5.10).

   Now by Definition 3.4, the number of colorings of
$G_e^*$ is $18 \cdot \psi(G^*, e)$.
For each coloring $\gamma^*$ of $G_e^*$, by
Theorem 3.3(C)(3), there are exactly $2 \cdot \psi(G',E)$
colorings $\gamma'$ of $G_E'$ such that the two equalities
in (5.10) hold.
Hence there are exactly $36 \cdot \psi(G^*,e), \psi(G',E)$
ordered pairs $(\gamma^*, \gamma')$ as in (5.10), and 
hence exactly that many colorings of $G_e$.  
Hence by Definition 3.4, eq.\ (5.1) holds.
That completes the proof of Theorem 5.3. 
\hfil ////\break

\noindent {\bf 6.\ \ Results of McKinney, and
results of Cappon and Walther} \hfil\break

   The paper [Br2], giving 
Theorem 5.1 and its proof, was published in March 2006.
A few months later, in the summer of 2006, in an
eight-week Mathematics 
REU (Research Experience for Undergraduates) program at Indiana University (organized by Professor Victor Goodman 
of the Indiana University Mathematics Department), 
Scott A.\ McKinney, at that time an undergraduate 
mathematics major at Cornell University, did some 
research on snarks (connected with the material in
this survey paper).
McKinney [McK] wrote a paper on the results of his 
research, as part of a collection of papers
by the students in that REU program.
His results were of a spirit similar to Theorem 5.1
and Theorem 5.3, and included an extension of 
Theorem 5.1 itself.  (More on that below.)
\hfil\break

   Seven years later, in the summer of 2013, in a 
similar eight-week Mathematics REU program at Indiana University
(organized by Professor Kevin Pilgrim of the Indiana
University Mathematics Department), further research
on snarks (again connected with the material in this
survey paper) was done jointly by two (then) 
undergraduate mathematics majors:
Ariana Cappon of Indiana University and Emily Walther of
Westminster College (of New Wilmington, Pennsylvania).  
Cappon and Walther [CpWl] wrote a joint paper on the
results of their research, as part of a
collection of papers by the students
in that REU program.  (More on that below.)
\hfil\break

   The REU research of all three students alluded to
above was mentored by
the author of this survey paper.
The research of Cappon and Walther [CpWl] built
directly on the earlier work of McKinney [McK], with 
strong encouragement from McKinney himself.
Here in Section 6, we shall summarize the results of
both papers.
We first describe the work of McKinney [McK], starting
with statements of the two main results of his paper:
\hfil\break 

   {\bf Theorem 6.1} (McKinney [McK, Theorem 5.5]).\ \
{\sl Suppose $j$ and $k$ are each a nonnegative integer.
Then there exists a cyclically 5-edge-connected snark $G$
and an edge $e$ of $G$ such that
$\psi(G,e) = 5^j \cdot 7^k$.} 
\hfil\break

   {\bf Theorem 6.2} (McKinney [McK, Corollary 5.6]).\ \
{\sl Suppose $j$, $k$, $\ell$, and $m$ are each a
nonnegative integer.
Then there exists a snark $G$ and an edge $e$ of $G$
such that
$\psi(G,e) = 2^j \cdot 3^k \cdot 5^\ell \cdot 7^m$.}
\hfil\break

   (In Theorem 6.2, the snark $G$ need not be cyclically
5-edge-connected, though by Definition 2.9 it must 
of course be at least cyclically 4-edge-connected.) 
\hfil\break

   Obviously Theorem 6.2 generalizes Theorem 5.1.
McKinney [McK] obtained Theorem 6.1 first, and then
mimicked the induction argument in [Br2] for Theorem 5.1
in order to derive Theorem 6.2 as a corollary of
Theorem 6.1.
Theorem 6.1 itself was obtained in [McK] by an induction argument that started with the Petersen graph and 
Theorem 3.5 and then (in the induction step) involved a particular choice of Kochol [Ko1] ``superposition''.
In Remark 6.3 below, we shall give just a brief
description of that whole process.
\hfil\break

   {\bf Remark 6.3.}\ \ (a) The context studied by 
McKinney [McK] was a follows:
It started with an arbitrary snark $G_0$.  
Two adjacent edges of that snark were ``replaced''
together by two (disjoint) Petersen graphs in a certain 
way, to create a new, ``bigger'' snark $G$.  
In that ``double superposition'' process (as we shall 
call it here for convenience), in each of those two
Petersen graphs, one edge was ``split'' into two ``strands''
and one vertex (not connected to that edge or
even adjacent to an end-point vertex of it) was ``split''
into three ``strands''; and the resulting ten ``strands'' 
were then ``tied up'' to each other and to (what was 
left of) $G_0$ in a particular way, creating a ``five edge
connection'' between (what was left of) $G_0$ and 
the union of (what was left of) the two Petersen graphs.  
Although the details were quite different, the general
spirit was somewhat like that of Context 5.2 (where one
edge of a snark was ``replaced'' by another snark,
in a ``superposition'' process in which two non-adjacent vertices of that ``other'' snark were each ``split'' into
three ``strands''). 
\medskip

   (b) McKinney [McK, Theorem 5.3] showed that for any 
given edge $e$ of the original snark $G_0$ that is not involved in the ``double superposition by two Petersen graphs'' (and hence $e$ is also an edge of the new 
snark $G$), one has that 
$\psi(G,e) = 5 \cdot \psi(G_0,e)$. 
\medskip

   (c) McKinney also studied the edge ${\bf e}$ of
$G_0$ that was (different from and) adjacent to the two
edges of $G_0$ that were ``replaced'' by Petersen graphs 
in the ``double superposition''.
In that ``double superposition'' process, the edge
${\bf e}$ itself was removed, and in the new snark $G$
(regardless of the original choice of $G_0$) it had  
a ``natural counterpart'' --- a new edge $E$.
McKinney [McK, Theorem 5.1] showed that
$\psi(G,E) = 7 \cdot \psi(G_0, {\bf e})$. 
\medskip

   (d) In proving both of the results described in
paragraphs (b) and (c) above, McKinney [McK] employed
arguments that were, while quite different in their
details, somewhat of the general spirit of the proof
of Theorem 5.3.
As mentioned above, McKinney [McK] then used induction,
starting with the Petersen graph and Theorem 3.5 and
then employing the results in both 
paragraphs (b) and (c) above in the induction 
step, to prove Theorem 6.1; and he then mimicked the
induction argument in [Br2] in order to derive
Theorem 6.2 as a corollary.
\hfil\break

   Now the rest of Section 6 here will be devoted to
a brief description of the work of Cappon and 
Walther [CpWl].
We start with a technical definition from their paper.
\hfil\break

   {\bf Definition 6.4.}\ \ Let $S_1$ denote the set 
of all prime numbers $p$ such that $p \leq 149$. 
Let $S_2$ denote the following set of prime 
numbers (greater than 149):
$S_2 := \{173, 179, 181, 197, \allowbreak 
229, 257, 271, 359\}$.
Define the set $S$ of prime numbers by 
$S := S_1 \cup S_2$.
\hfil\break

   Cappon and Walther [CpWl] started with the scheme 
studied by McKinney, systematically tried certain
sequences of embellishments (extra vertices and edges), and devised and employed an efficient algorithm for recursively keeping a running track of the numbers of possible 
edge-3-colorings resulting therefrom.
Thereby they derived a version of the result of McKinney
described in Section 6.3(b) above, but with the factor
5 there replaced by arbitrary greater prime numbers in
the set $S$ above.
Using induction, building on the arguments in
[Br2] and [McK] for Theorem 5.1 and Theorem 6.2, 
they extended those theorems to the following form:
\hfil\break

   {\bf Theorem 6.5} (Cappon and Walther
[CpWl, Theorem 4.2]).
{\sl Suppose that for each $p \in S$ (from 
Definition 6.4 above), $m(p)$ is a nonnegative integer.
Then there exist a snark $G$ and en edge $e$ of $G$
such that $\psi(G,e) = \prod_{p \in S} p^{m(p)}$.}
\bigskip

   {\bf Remark 6.6.} (a) It was clear that the
embellishments of the scheme of McKinney studied by 
Cappon and Walther [CpWl] could have led to 
many more prime numbers in Theorem 6.5 besides the 
ones in the set $S$; but time ran out in their 
eight-week REU
program, and the report [CpWl] had to be written.
\medskip

   (b) Refer to Definition 4.6.
Cappon and Walther [CpWl, Theorem 5.2] also
proved that under the hypothesis (first sentence)
of Theorem 6.5, along with the extra assumption
$m(3) \neq 1$, there exist a snark $G$ and a 
pentagon $P$ in $G$ such that     
$\psi(G,P) = \prod_{p \in S} p^{m(p)}$.
In that result involving pentagons,
the number 1 (the case $m(p) = 0$ for all $p \in S$)
comes from Theorem 3.5,  
the powers of 2 come from a simple observation in the 
context of the proof of [Br2, Theorem 2.2],
and the powers of $5$ come from a simple 
observation in the context of an argument of McKinney 
[McK, Theorem 5.3] (see Remark 6.3(b) above).
The recursive algorithm of Cappon and Walther [CpWl]
alluded to in the first sentence after Definition 6.4 
above, took care (for pentagons) of powers of all prime numbers $p \in S$ 
such that $p \geq 7$, and also 
(via an extra induction argument) took care 
of powers of 3 with exponent $m(p) \geq 2$ (and 
of course the case $m(3) = 0$ is covered trivially) --- 
but not $m(3) = 1$.
In particular, it is unknown whether there exists a 
snark $G$ with a pentagon $P$ such that $\psi(G,P) = 3$.
\medskip

   (c) Cappon and Walther [CpWl] also studied a variation
of McKinney's scheme in which, for an arbitrary snark $G_0$,
three adjacent edges, all connected to the same vertex,
are replaced together by three Petersen graphs, to
created a ``bigger'' snark $G$.
With an argument similar to that of McKinney described 
in Remark 6.3(b) above, 
Cappon and Walther [CpWl, Section 5.2 
(see final two paragraphs)]
showed that for any 
given edge $e$ of the original snark $G_0$ that is not involved in the ``triple superposition by three Petersen graphs'' (and hence $e$ is also an edge of the new 
snark $G$), one has that 
$\psi(G,e) = 11 \cdot \psi(G_0,e)$.
(With an extra embellishment --- two edges inserted --- 
they also derived in this scheme the same result with the factor 11 replaced by 19.)
\medskip

   (d) With yet another variation on McKinney's scheme,
Cappon and Walther [CpWl, Theorem 5.4] showed that there
are infinitely many prime numbers $p$ such that $p$ is a
divisor of $\psi(G,e)$ for some snark $G$ and some edge 
$e$ of $G$.  
(That infinite set of prime numbers was not determined
explicitly; the argument was ``nonconstructive''.)

\bigskip

\noindent {\bf 7.\ \ Some open problems} \hfil\break

   {\bf Section 7.1.  Some open problems.}\ \ Given 
below is a list of some open problems involving snarks 
and the numbers $\psi(G,e)$, and also involving 
edge-3-colorable cubic graphs with orthogonal edges.  
These problems are motivated by the papers of 
K\'aszonyi [K\'a1, K\'a2, K\'a3] and are rooted primarily
in the material in Section 3 of this survey paper.
The first problem was implicitly posed by the
author in [Br2] and was mentioned in Section 5 above. \hfil\break

   {\bf Problem 1.}\ \ For precisely what positive 
integers $n$ do there exist a snark $G$ and an edge 
$e$ of $G$ such that $\psi(G,e) = n$? \hfil\break

   This problem can perhaps be approached via the following closely related one: \hfil\break

   {\bf Problem 2.}\ \ For what prime numbers $p$ does the
following (uncertain) ``Hypothesis $S(p)$'' hold? \smallskip

\noindent {\it Hypothesis $S(p)$}: If $n$ is a positive integer, $G$ is a snark, $e$ is an edge of $G$, and 
$\psi(G,e) = n$, then there exist a snark ${\cal G}$ and an edge $E$ of ${\cal G}$ such that 
$\psi({\cal G},E) = n \cdot p$. \hfil\break

   Now Hypothesis $S(p)$ was (implicitly) verified for 
$p = 2, 3$ by [Br2, Theorem 2.2] (one combines that 
result with Theorem 3.5, using the Petersen graph); 
it was verified for $p = 5,7$ by 
McKinney [McK, Theorems 5.3 and 5.1]
(see Remark 6.3(b)(c) above), and it was verified
by Cappon and Walther [CpWl, Lemma 4.3] for all
higher prime numbers $p$ in the set $S$ in 
Definition 6.4.   
The next two problems are variations on Problem 1, 
and are motivated by Theorem 6.1.  
Perhaps they can be approached 
by corresponding variations on Problem 2. 
(For Problem 3, recall the result of McKinney in 
Theorem 6.1.) \hfil\break 

   {\bf Problem 3.}\ \ For precisely what positive 
integers $n$ do there exist a cyclically 5-edge-connected snark $G$ and an edge $e$ of $G$ such that $\psi(G,e) = n$? \hfil\break  
              
   {\bf Problem 4.}\ \ For precisely what positive 
integers $n$ do there exist a cyclically 6-edge-connected snark $G$ and an edge $e$ of $G$ such that $\psi(G,e) = n$? \hfil\break

   Problem 4 is motivated partly by the ``flower'' snarks
in Remark 2.14.
The flower snarks $J_n$ for $n \in \{7,9,11,\dots\}$
(but not $J_5$) 
are known to be cyclically 6-edge-connected. 
The next question seems to be of obvious special 
interest: \hfil\break

   {\bf Problem 5.}\ \ What are the numbers $\psi(G,e)$ for the flower snarks $G$ and their edges $e$? \hfil\break

   For a given flower snark $G$, the edges fall into four
equivalence classes in which two edges $e_1$ and $e_2$ are
``equivalent'' if there is an automorphism of $G$ in which $e_1$ is mapped to $e_2$.
Of course $\psi(G,e_1) = \psi(G,e_2)$ for any two such ``equivalent'' edges $e_1$ and $e_2$.
(Recall the Remark after Definition 3.4.)\ \ 
This suggests that as a solution of Problem 5, there might 
be a ``recursion formula'' that involves as a parameter 
the subscript $n \in \{5,7,9,11,\dots\}$ 
for a given flower snark $J_n$
(in the notations of Remark 2.14), 
and that also involves the
four numbers $\psi(J_n,e)$ corresponding to the four ``equivalence classes'' of edges $e$ of $J_n$. \hfil\break

   The next two problems involve the material in Section 4.
\hfil\break

   {\bf Problem 6.}\ \ Refer to Definition 4.6.  For what positive integers $n$ do there exist a snark $G$ and a pentagon $P$ in $G$
such that $\psi(G,P) = n$? \hfil\break

   The progress on this problem so far is summarized in
Remark 6.6(b). 
\bigskip

   {\bf Problem 7.}\ \ In Theorem 4.8, five edges of $G$ 
were omitted: the five ``connecting'' edges 
$(t_k, w_{2k})$, $k \in {\bf Z}_5$.
If $e$ is one of those five edges, what can one 
say about $\psi(G,e)$ in terms of $\psi(G',e')$ and 
$\psi(G^*, e^*)$ for appropriate edges $e'$ and $e^*$ of (respectively) $G'$ and $G^*$? \hfil\break

   From the solutions to Problems 5 and 7, or perhaps
more easily from a direct argument, one might compute the numbers $\psi(G,e)$ for all of the edges $e$ of  
Isaac's [Is] ``double star'' snark $G$ alluded to right 
before Context 4.7.
\hfil\break

   The next problem involves the following definition:
A snark $G$ is said to be ``critical'' if
$\psi(G,e) \geq 1$ (i.e.\ $G_e$ can be edge-3-colored)
for every edge $e$ of $G$.
Not all snarks are critical.
(In fact some snarks $G$ are so severely ``anti-critical'' 
that $\psi(G,e) = 0$ for {\it every\/} edge of $G$; for further information and earlier references on such 
snarks, see [BGHM, Section 4.7] and [H\"a, Section 3].)\ \ 
Starting with the Petersen graph and Theorem 3.5, and
applying induction using [Br2, Theorem 2.2 and subsequent
sentence], one obtains Theorem 5.1 in Section 5 above 
with the word ``snark'' replaced by the phrase
``critical snark''.
That suggests the following problem:
\hfil\break   

   {\bf Problem 8.}\ \ For what positive integers $n$
do there exist a {\it critical\/} snark $G$ and an edge
$e$ of $G$ such that $\psi(G,e) = n$?
\hfil\break

   The remaining problems below come directly (at least
implicitly) from the work of K\'aszonyi [K\'a1, K\'a2].
\hfil\break   

   Refer to Definition 3.1.  In the proof of Theorem 3.5,
for the Petersen graph ${\cal P}$ and an edge $e$ of 
${\cal P}$, the graph ${\cal P}_e$ was represented (as in
[K\'a1, K\'a2]) as an ``8-vertex wheel with four 
rim-to-rim spokes''.
In the notations used there (in the proof of Theorem 3.5)
for that graph, the orthogonal edges resulting directly
from the ``removal'' of the edge $e$ from ${\cal P}$
were denoted $f_0$ and $f_2$.
(Those were the edges corresponding to
$d_1$ and $d_2$ in Notations 3.2.)\ \ 
However, by simple symmetry (simply ``rotate the wheel
45 degrees''), that graph ${\cal P}_e$ has another pair
of orthogonal edges: $f_1$ and $f_3$.
This suggests the following problem: \hfil\break

   {\bf Problem 9.}\ \ If $G$ is a snark, $e$ is an
edge of $G$, and $G_e$ can be colored 
(i.e.\ $\psi(G,e) \geq 1$), does the cubic graph $G_e$ 
have (at least) two pairs of orthogonal edges?
Or are there instead examples where $G_e$ has only one
pair of orthogonal edges (the pair  
identified by K\'aszonyi [K\'a2] in Theorem 3.3(B)
--- the edges $d_1$ and $d_2$ in Notations 3.2)?
If the latter is the case, then under what extra
assumptions on the snark $G$ and the edge $e$ of $G$ does
$G_e$ have at least two pairs of orthogonal edges?
(Just one pair?) \hfil\break     

   For the final question, Problem 10 below, a definition 
will be given first:
Suppose $G$ is a simple cubic graph which, say, is 
(at least) cyclically 4-edge-connected and has girth 
at least 5.  
(No assumption on whether or not $G$ can be 
edge-3-colored.)\ \  Suppose $e$ is an edge of $G$.  
Let $d_1$ and $d_2$ be the edges of $G_e$ specified in
Notations 3.2.  
Let us say that the edge $e$ satisfies Condition
${\cal K}$ (for K\'aszonyi) if 
(i) the (simple cubic) graph $G_e$
can be edge-3-colored, and 
(ii) the edges $d_1$ and $d_2$ are
orthogonal (again see Definition 3.1).  \hfil\break

   {\bf Problem 10.}\ \ Suppose $G$ is a simple cubic graph
which (say) is (at least) cyclically 4-edge-connected and
has girth at least 5.
If some edge of $G$ satisfies Condition ${\cal K}$ (see the
preceding paragraph above), does it follow that $G$ is a
snark? \hfil\break

   This question is quite specific.
If the answer is ``no'', then there are obvious
variations on this question.
For example, what if at least two edges of $G$ satisfy
Condition ${\cal K}$?   
If the answer is still ``no'', then (for example) 
what if all five edges of some pentagon (if one exists)
in $G$ satisfy Condition ${\cal K}$? 
\hfil\break 
   
{\bf Section 7.2. Final Remarks.}\ \ Here are some final comments on the papers of 
K\'aszonyi [K\'a1, K\'a2, K\'a3] on which this survey 
paper is based.
\medskip

   (A) For a long time, those three papers of K\'aszonyi 
did not seem to be known much in the ``snark community''. 
(The author of this survey paper has not found
any citations to those papers of K\'aszonyi
in other published papers
prior to their citations in the 2006 paper [Br2].) 
\medskip
  
   (B) In Section 5 (its first paragraph and a couple of
other places), the paper of Isaacs [Is] is cited for a
``dot product'' of two snarks ---  a particular 
``4-edge-connection'' procedure for combining two 
``smaller'' snarks to form a ``bigger'' one.  
In fact a few years earlier,
K\'aszonyi [K\'a1, pp.\ 86-87, Operation 3]
had presented an exactly analogous
``4-edge-connection'' for combining two simple cubic
graphs, each of them being edge-3-colorable with
a pair of orthogonal edges (one of those graphs being
${\cal P}_e$ for an edge $e$ of the Petersen 
graph ${\cal P}$), to form a ``bigger'' simple
cubic graph which is edge-3-colorable with a pair of
orthogonal edges. 
\medskip

   (C) {\bf Acknowledgement of priority.}\ \ 
Certain results and arguments of K\'aszonyi [K\'a2, K\'a3] 
--- roughly (recall the first paragraph of Section 3),   
statements and proofs of 
Theorem 3.3(B), Theorem 3.3(C)(1)(2), and 
Theorem 3.7(B) --- were independently rediscovered 
a few years later by the 
author [Br1, Theorem 1, Lemmas 1, 2, and 8,  
Corollary 2, and their proofs]. 
The priority for those results and arguments 
belongs to K\'aszonyi. 
\medskip

   (D) The paper [K\'a3] is somewhat cryptic.  For a given 
edge-3-colorable simple cubic graph $H$ with orthogonal 
edges, K\'aszonyi [K\'a3] defined a ``coloring graph'', 
which will be referred to here as ${\cal H}$.
The ``vertices'' of ${\cal H}$ correspond to 
edge-3-colorings of $H$ with three given ``colors'' (say 
the elements $a$, $b$, and $c$ from (2.1)).  
Two ``vertices'' of ${\cal H}$ are connected by an ``edge''
of ${\cal H}$ if the two corresponding edge-3-colorings of 
$H$ differ from each other by just the interchanging of the two colors on a single Kempe cycle.
The observations made by K\'aszonyi [K\'a3, p.\ 35] that (cryptically) yielded Theorem 3.3(A)(C) and Theorem 3.7 
were made in the terminology of ``coloring graphs''.  
In giving those arguments of K\'aszonyi here (in the proofs 
of Theorems 3.3 and 3.7), we have simply transcribed  
K\'aszonyi's own presentation of those arguments, involving
the terminology of ``coloring graphs'', into the more transparent terminology of edge-3-colorings and Kempe 
cycles. 
\medskip

   (E) To summarize, the work of 
K\'aszonyi [K\'a1, K\'a2, K\'a3], along with some of the related later work of other people as described above, 
provide a collection of mathematical problems
(including, but not limited to, the ones listed above)
which can be attacked without too much specialized mathematical preparation, and which are in particular 
well suited for independent research projects for undergraduate mathematics students.
This survey paper can hopefully
facilitate research on such problems. 
\hfil\break

   {\bf Acknowledgement.}\ \ The author thanks 
Scott McKinney for his proofreading of this manuscript, 
for his valuable suggestions which helped improve the
exposition, and for calling attention to recent
pertinent work on snarks such as in the references [BGHM]
and [H\"a].
\hfil\break 

\centerline {\bf References} \bigskip  

\itemitemitem {[Br1]}  R.C.\ Bradley.  A remark on noncolorable
cubic graphs.  {\it J.\ Combin.\ Theory Ser B\/} 24
(1978) 311-317. \medskip

\itemitemitem {[Br2]}  R.C.\ Bradley.  On the number of colorings 
of a snark minus an edge.  {\it J.\ Graph Theory\/} 51
(2006) 251-259.  \medskip

\itemitemitem {[BGHM]}  G.\ Brinkmann, J.\ Goedgebeur,
J.\ H\"agglund, and K.\ Markstr\"om.
Generation and properties of snarks.
arXiv:1206.6690v2 [math.CO] 5 Nov 2012.  \medskip

\itemitemitem {[CCW]}  P.J.\ Cameron, A.G.\ Chetwynd, and
J.J.\ Watkins.  Decompositions of snarks.
{\it J.\ Graph Theory\/} 11 (1987) 13-19.  \medskip

\itemitemitem {[CpWl]}  A.\ Cappon and E.\ Walther.
Prime factorization of K\'aszonyi numbers.
In: {\it Research Experience for Undergraduates,
Research Reports, Indiana University,
Bloomington, Summer 2013\/}, pp.\ 60-89.
Indiana University, Bloomington, Indiana, Summer 2013.
Posted online at \hfil\break
www.math.indiana.edu/reu/2013/reu2013.pdf 
\medskip

\itemitemitem {[Ga]}  M. Gardner.  Mathematical games:  Snarks,
boojums, and other conjectures related to the four-color
map theorem.  {\it Sci.\ Amer.\/} 234 (1976) 126-130.
\medskip

\itemitemitem {[H\"a]}  J.\ H\"agglund.
On snarks that are far from being 3-edge-colorable.
\hfil\break
arXiv:1203.2015v1 [math.CO] 9 Mar 2012.  \medskip

\itemitemitem {[Is]}  R.\ Isaacs.  Infinite families of nontrivial
trivalent graphs which are not Tait colorable.
{\it Amer.\ Math.\ Monthly\/} 82 (1975) 221-239.  \medskip

\itemitemitem {[K\'a1]}  L.\ K\'aszonyi.  A construction of cubic
graphs, containing orthogonal edges.
{\it Ann.\ Univ.\ Sci.\ Budapest E\"otv\"os Sect.\ Math.\/}
15 (1972) 81-87. \medskip

\itemitemitem {[K\'a2]}  L.\ K\'aszonyi.  On the 
nonplanarity of some cubic graphs.
{\it Ann.\ Univ.\ Sci.\ Budapest E\"otv\"os Sect.\ Math.\/}
15 (1972) 123-131. \medskip

\itemitemitem {[K\'a3]}  L.\ K\'aszonyi.  On the structure 
of coloring graphs.
{\it Ann.\ Univ.\ Sci.\ Budapest E\"otv\"os Sect.\ Math.\/}
16 (1973) 25-36. \medskip

\itemitemitem {[Ko1]}  M.\ Kochol.  Snarks without small cycles.
{\it J.\ Combin.\ Theory Ser B\/} 67 (1996) 34-47. 
\medskip

\itemitemitem {[Ko2]}  M.\ Kochol.  Superposition and
construction of graphs without nowhere-zero $k$-flows.
{\it Eur.\ J.\ Combin.\/} 23 (2002) 281-306.  \medskip

\itemitemitem {[McK]} S.A.\ McKinney.  On the number of 
edge-3-colourings of a snipped snark.
In: {\it Research Experiences for Undergraduates,
Student Reports, Indiana University\/}, pp.\ G1-G15.
Indiana University, Bloomington, Indiana, Summer 2006.
(Reference copy available at the Swain Hall Library,
Indiana University, Bloomington, Indiana.)\ \ 
Posted in slightly embellished form on \hfil\break
arXiv:1304.5427v1 [math.CO] 19 Apr 2013 \medskip

\itemitemitem {[TW]} F.C.\ Tinsley and J.J.\ Watkins.
A study of snark embeddings.
{\it Graphs and Applications (Boulder, Colorado, 1982)\/},
pp.\ 317-332.
Wiley, New York, 1985. \medskip

\itemitemitem {[Wi]}  R.\ Wilson.  {\it Four Colors 
Suffice.\/}
Princeton University Press, Princeton, 2002. \medskip 

\bye